\documentclass[10pt]{amsart}
\usepackage[utf8]{inputenc}
\usepackage[T1]{fontenc}
\usepackage[english]{babel}
\usepackage{color}

\usepackage{comment}
\usepackage{amsmath}
\usepackage{amssymb}
\usepackage{dsfont}
\usepackage{amsfonts}
\usepackage{amsthm}
\usepackage{tikz}
\usepackage{graphicx}
\usepackage{fancyhdr}
\usepackage{caption}
\usepackage{enumerate}
\usepackage{mathtools}

\usepackage[colorlinks=true, pdfstartview=FitH, bookmarks=false]{hyperref}
\hypersetup{urlcolor=black, linkcolor=black, citecolor=black, colorlinks=true}

\newcommand{\vertiii}[1]{{\left\vert\kern-0.25ex\left\vert\kern-0.25ex\left\vert #1 https://www.overleaf.com/project/57d27f441ebc8db50b18ed7d
    \right\vert\kern-0.25ex\right\vert\kern-0.25ex\right\vert}}

\newcommand\Z{\mathbb{Z}}

\theoremstyle{plain}
\newtheorem{definition}{Definition}[section]
\newtheorem{proposition}[definition]{Proposition}
\newtheorem{corollary}[definition]{Corollary}
\newtheorem{theorem}[definition]{Theorem}
\newtheorem{lemma}[definition]{Lemma}
\newtheorem*{thm*}{Theorem}
\newtheorem*{prop*}{Proposition}
\newtheorem*{lem*}{Lemma}

\theoremstyle{remark}
\newtheorem{remark}[definition]{Remark}

\newtheorem*{rem*}{Remark}

\numberwithin{equation}{section}

\usepackage{etoolbox}
\apptocmd{\sloppy}{\hbadness 10000\relax}{}{}
\apptocmd{\sloppy}{\vbadness 10000\relax}{}{}

\title[LLT on relatively hyperbolic groups wrt virt. nilpotent subgroups]{Local limit theorems on relatively hyperbolic groups with respect to virtually nilpotent subgroups}
\author{Matthieu Dussaule}

\date{}

\begin{document}

\begin{abstract}
Given a probability measure on a finitely generated group, the local limit problem consists in finding asymptotics of $p_n(e,e)$, the probability that the random walk at time $n$ is at the origin.
We give the classification of all possible local limit theorems, up to bounded error, for finitely supported, symmetric, admissible probability measures on a relatively hyperbolic group with respect to virtually nilpotent subgroups.
\end{abstract}

\maketitle

\section{Introduction}

\subsection{Random walks and the local limit problem}
Let $\Gamma$ be a finitely generated group and let $\mu$ be a probability measure on $\Gamma$.
The random walk driven by $\mu$ and starting at the neutral element $e$ is the sequence of random variables $X_n$, which is defined by
$$X_n=g_1...g_n,$$
where the $g_k$ are independent random variables whose common distribution is given by $\mu$.
The random walk is called symmetric if $\mu(g)=\mu(g^{-1})$ for all $g\in \Gamma$.
It is called admissible if the support of $\mu$ generates $\Gamma$ as a semi-group.

The $n$th convolution power $\mu^{(n)}$ of $\mu$ is defined by $\mu^{(0)}(x)=\delta_e(x)$, $\mu^{(1)}(x)=\mu(x)$ and for $n\geq 2$
$$\mu^{(n)}(x)=\sum_{x_1,...,x_{n-1}}\mu(x_1)\mu(x_1^{-1}x_2)\mu(x_2^{-1}x_3)...\mu(x_{n-1}^{-1}x).$$
The probability measure $\mu^{(n)}$ is the distribution of $X_n$, the random walk at time $n$.
We will sometimes write $p_n(x,y)=\mu^{(n)}(x^{-1}y)$.
An admissible random walk is called aperiodic if there exists $n_0$ such that for every $n\geq n_0$, $\mu^{(n)}(e)>0$.

We then introduce the Green function, which is the power series defined by
$$G(x,y|r)=\sum_{n\geq 0}p_n(x,y)r^n.$$
As we will see below, this function encodes many asymptotics properties of the random walk.
We denote by $R$ its radius of convergence.
Its inverse $\rho=R^{-1}$ is called the spectral radius of the random walk and it satisfies $\rho=\limsup \mu^{(n)}(e)^{1/n}$.

\medskip
Given a random walk on a finitely generated group, the local limit problem consists in finding asymptotics of
$p_n(x,y)$, as $n$ goes to infinity.
In several standard examples, one has a local limit theorem of the form
\begin{equation}\label{e:standardllt}
p_n(e,x)\sim C_xR^{-n}n^{-\alpha},
\end{equation}
where $C_x$ is a constant depending only on $x$, $R$ is the inverse of the spectral radius of the random walk and $\alpha$ is an exponent sometimes called the critical exponent of the random walk.

For example, given an admissible probability measure with finite exponential moments on a finitely generated abelian group of rank $d$,~(\ref{e:standardllt}) holds with $\alpha=d/2$, see \cite[Theorem~13.12]{Woess-book}.
Also, given a finitely supported, admissible, symmetric probability measure on a hyperbolic group,~(\ref{e:standardllt}) holds with $\alpha=3/2$, see \cite{GouezelLLT} and references therein.

\medskip
In \cite{Gerl}, P.~Gerl conjectured that if~(\ref{e:standardllt}) holds for a finitely supported random walk, then $\alpha$ is a group invariant.
This conjecture was disproved by
Cartwright in \cite{Cartwright88}, \cite{Cartwright89}.
He gave examples of local limit theorems on free products $\Z^d*\Z^d$, with $\alpha=d/2$ and examples on the same groups with $\alpha=3/2$.

Free products are, in fact, combinatorial models for relatively hyperbolic groups, whose precise definition is recalled in Section~\ref{ss:defrh} and these groups are a great source of examples of particular behaviors for random walks.
Recently, Cartwright's results have been generalized in a series of papers \cite{DussauleLLT1}, \cite{DussauleLLT2}, \cite{DPT}, \cite{DPT23}, which deal with the local limit problem on relatively hyperbolic groups.
The main idea established in these papers is that there is a competition between the underlying hyperbolic structure and the so-called parabolic subgroups, which can be read in the local limit theorem.

Our main goal in this paper is to end the classification of local limit theorems, up to bounded error, on relatively hyperbolic groups whenever parabolic subgroups are virtually nilpotent.

\medskip
Given two positive functions $f$ and $g$, we will write $f\asymp g$ if the ratio of $f$ and $g$ is bounded from above and from below, that is, if there exists $C$ such that $C^{-1}f\leq g \leq C f$.
We now state the following theorem, which summarizes the classification.
All relevant terminology will be introduced below.

\begin{theorem}\label{t:maintheoremsummary}
Let $\Gamma$ be a non-elementary relatively hyperbolic group with respect to virtually nilpotent subgroups.
Let $\mu$ be a finitely supported, admissible and symmetric probability measure on $\Gamma$.
Assume that $\mu$ is aperiodic.
Then the following holds.
\begin{itemize}
\item If $\mu$ is convergent, then
    $$p_n(e,e)\asymp R^{-n}n^{-d/2},$$
    where $d$ is the rank of spectral degeneracy of $\mu$.
    \item If $\mu$ is divergent and spectrally positive recurrent, then
    $$p_n(e,e)\asymp R^{-n}n^{-3/2}.$$
    \item If $\mu$ is divergent and not spectrally positive recurrent, then its rank of spectral degeneracy $d$ must be 5 of 6.
    \begin{itemize}
        \item If $d=5$, then
        $$p_n(e,e)\asymp R^{-n}n^{-5/3}.$$
        \item If $d=6$, then
        $$p_n(e,e)\asymp R^{-n}n^{-3/2}\big(\log (n)\big)^{-1/2}.$$
    \end{itemize}
\end{itemize}
If $\mu$ is not aperiodic, then similar asymptotics hold for $p_{2n}(e,e)$.
\end{theorem}

Roughly speaking, in the convergent case, only the parabolic subgroups are influential and in the spectrally positive recurrent case, only the underlying hyperbolic structure is influential.
In the third case, both are influential.
Before we give more explanations about this classification, let us give a short overview of previous results from \cite{DussauleLLT1}, \cite{DussauleLLT2}, \cite{DPT}, \cite{DPT23}.

In \cite{DussauleLLT1}, it was proved that in the spectrally positive recurrent case, one has $p_n(e,e)\asymp R^{-n}n^{-3/2}$.
It turns out that one does not need to assume that parabolic subgroups are virtually nilpotent.
In \cite{DussauleLLT2}, this result was improved to the precise asymptotics $p_n(e,x)\sim C_xR^{-n}n^{-3/2}$, assuming further that $\mu$ is not spectrally degenerate. We will discuss this property below, let us only mention here that spectral non-degeneracy implies spectral positive recurrence.

In \cite{DPT}, it was proved that in the convergent case $p_n(e,x)\sim C_xR^{-n}n^{-d/2}$, assuming that parabolic subgroups are virtually abelian.
The situation with virtually nilpotent subgroups was briefly discussed, but the first item of Theorem~\ref{t:maintheoremsummary} is new.

It was also noted in \cite{DPT} that the classification of local limit theorems on free products, as it was known so far from \cite{CandelleroGilch}, was incomplete and that the case of a divergent, not spectrally positive recurrent probability measure could occur.
This led to the paper \cite{DPT23} in which such probability measures are constructed and where the asymptotics in the third item of Theorem~\ref{t:maintheoremsummary} are proved on one specific example.
Proving here that these asymptotics still hold in general is also new.

In other words, the main contributions of this paper are the first and third items of Theorem~\ref{t:maintheoremsummary}.

\subsection{The fundamental asymptotic differential equation}\label{ss:fundamentalasymptoticequation}
We consider a relatively hyperbolic group $\Gamma$ endowed with a probability measure $\mu$.
We introduce the following quantities :
$$I^{(1)}(r)=\sum_{x\in \Gamma}G(e,x|r)G(x,e|r)$$
and
$$I^{(2)}(r)=\sum_{x,y\in \Gamma}G(e,x|r)G(x,y|r)G(y,e|r).$$

Also given a parabolic subgroup $H$ we set
$$I_H^{(2)}(r)=\sum_{x,y\in H}G(e,x|r)G(x,y|r)G(y,e|r).$$
Choosing a full finite set $\Omega_0$ of representatives of conjugacy classes of parabolic subgroups, we then define
$$J^{(2)}(r)=\sum_{H\in \Omega_0}I^{(2)}_H(r).$$
As we will see, these quantities are intimately related to the derivatives of the Green function.
In particular~(\ref{e:I1G1}) and~(\ref{e:I2G2}) below show that
$$I^{(1)}(r)\asymp G'(e,e|r)$$
and
$$I^{(2)}(r)\asymp G''(e,e|r).$$

In many situations, local limit theorems are derived from the study of singularities of the Green function.
Typically, whenever it is possible to find a precise expansion of the Green function with error terms, one may apply Darboux theory to deduce the precise asymptotics of $\mu^{(n)}(e)$.
This method has been used a lot to prove local limit theorems on free groups and free products and is presented in detail in \cite[Chapter~17, Chapter~19]{Woess-book}.
In \cite{GouezelLalley} and \cite{GouezelLLT}, the authors developed a very useful method based on Karamata's Tauberian Theorem to find asymptotics of $\mu^{(n)}(e)$ in hyperbolic groups using only asymptotics of $G(e,e|r)$ and its derivatives as the real variable $r$ tends to $R$.
The main strategy of \cite{GouezelLalley} is then to find an asymptotic differential equation satisfied by the Green function to deduce that $G'(e,e|r)\sim C\frac{1}{\sqrt{R-r}}$.
Namely, it is proved that $I^{(2)}(r)\sim C I^{(1)}(r)^3$ or equivalently $G''(e,e|r)\sim C G'(e,e|r)^3$ as $r$ tends to $R$.

This method has been developed further on relatively hyperbolic groups.
The main result of \cite{DussauleLLT1} is that
\begin{equation}\label{e:fundamental}
I^{(2)}(r)\asymp J^{(2)}(r)I^{(1)}(r)^3.
\end{equation}
We call this the fundamental asymptotic differential equation.
All the local limit theorems proved in \cite{DussauleLLT1}, \cite{DussauleLLT2}, \cite{DPT} and \cite{DPT23} are derived from it.

\subsection{Main results}
We now give more details on the terms \emph{convergent}, \emph{divergent} and \emph{spectrally positive recurrent} in Theorem~\ref{t:maintheoremsummary}. We also explain how the main results of the paper are obtained.

In view of~(\ref{e:fundamental}), the two quantities $I^{(1)}(r)$ and $J^{(2)}(r)$ are of special importance and encode many properties of random walks on relatively hyperbolic groups.
\begin{definition}
Let $\Gamma$ be a non-elementary relatively hyperbolic group and let $\mu$ be a probability measure on $\Gamma$.
\begin{itemize}
    \item The random walk, or equivalently the probability measure $\mu$, is called divergent if $I^{(1)}(R)$ is infinite and convergent otherwise.
    \item The random walk, or equivalently the probability measure $\mu$, is said to have finite Green moments if $J^{(2)}(R)$ is finite.
    \item The random walk, or equivalently the probability measure $\mu$, is called spectrally positive recurrent if it is divergent and has finite Green moment.
\end{itemize}
\end{definition}

This terminology was introduced in \cite{DussauleLLT1} and comes from the close analogy between local limit theorems on relatively hyperbolic groups on the one hand and the study of the geodesic flow on negatively curved manifolds on the other hand.
We refer to \cite[Section~1.4]{DPT} and \cite[Section~3.3]{DussauleLLT1} for more details about this analogy.

We will also introduce below the notion of spectral degeneracy along a parabolic subgroup.
Roughly speaking, a probability measure $\mu$ on a relatively hyperbolic group $\Gamma$ is called spectrally degenerate along a parabolic subgroup $H$ if the first return kernel to $H$ reaches its spectral radius when $\mu$ reaches its own spectral radius; see Definition~\ref{d:spectraldegeneracy} for a precise definition.

It turns out that if $\mu$ is not spectrally degenerate along any parabolic subgroup, then it is spectrally positive recurrent.
In particular, if $\mu$ is convergent, it must be spectrally degenerate along one parabolic subgroup.
Whenever parabolic subgroups are virtually nilpotent, they all have a well-defined homogeneous dimension (see Section~\ref{ss:nilpotentgeometry} for more details).

In such a situation, we define the rank of spectral degeneracy $d$ of $\mu$ as the smallest possible homogeneous dimension of a virtually nilpotent parabolic subgroup along which $\mu$ is spectrally degenerate.
Note that this definition also makes sense if $\mu$ is not spectrally degenerate or if all parabolic subgroups are not virtually nilpotent, up to setting $d=+\infty$.

The main results of this paper are the following two theorems which restate, respectively, the first and third items of the classification given in Theorem~\ref{t:maintheoremsummary}.

\begin{theorem}\label{t:maintheoremconvergent}
Let $\Gamma$ be a non-elementary relatively hyperbolic group with respect to virtually nilpotent subgroups and let $\mu$ be a finitely supported symmetric admissible probability measure on $\Gamma$.
Assume that $\mu$ is convergent and aperiodic.
Let $d$ be the rank of spectral degeneracy of $\mu$.
Then,
$$p_n(e,e)\asymp R^{-n}n^{-d/2}.$$
If $\mu$ is not aperiodic, similar asymptotics hold for $p_{2n}(e,e)$.
\end{theorem}

\begin{theorem}\label{t:maintheoremcritical}
Let $\Gamma$ be a non-elementary relatively hyperbolic group with respect to amenable subgroups and let $\mu$ be a finitely supported symmetric admissible probability measure on $\Gamma$.
Assume that $\mu$ is divergent and has infinite Green moments and that it is aperiodic.
Let $d$ be the rank of spectral degeneracy of $\mu$.
Then, $d=5$ or $d=6$.
Moreover, if $d=5$, then
$$p_n(e,e)\asymp R^{-n}n^{-5/3}$$
and if $d=6$, then
$$p_n(e,e)\asymp R^{-n}n^{-3/2}\big(\log (n)\big)^{-1/2}.$$
If $\mu$ is not aperiodic, similar asymptotics hold for $p_{2n}(e,e)$.
\end{theorem}

\begin{remark}
It turns out one does not need to assume \textit{a priori} that some of the parabolic subgroups are virtually nilpotent for this second result.
Indeed, this is a consequence of the fact that $J^{(2)}(R)$ is infinite, see Corollary~\ref{c:J2infiniteimpliesnilpotency} for more details.
\end{remark}

Let us outline the proofs of these two theorems.
In the convergent case (Theorem~\ref{t:maintheoremconvergent}), the fundamental asymptotic differential equation yields
$I^{(2)}(r)\asymp J^{(2)}(r)$.
If $J^{(2)}(R)$ is finite, then this only says that $I^{(2)}(R)$ is also finite.
However, if $J^{(2)}(R)$ is infinite, we deduce that $I^{(2)}(R)$ is infinite and that $I^{(2)}(r)$ has the same behavior as $J^{(2)}(r)$.
With the interpretation in terms of derivatives of Green functions mentioned above, we can restate this as follows.
The second derivative of the Green function associated with $\mu$ is infinite at the spectral radius. Moreover, its behavior is controlled by the behavior of the second derivatives of induced Green functions on parabolic subgroups.
We refer to Section~\ref{ss:derivativesGreenfunction} for more details on the connection of the sums $I^{(k)}(r)$ and $J^{(k)}(r)$ with the derivatives of suited Green functions.
As we will see, the induced Green functions in question are, in fact, associated with the first return kernels to parabolic subgroups.
Now the quantity $J^{(2)}(R)$ might be finite, but using that parabolic subgroups are virtually nilpotent, it can be proved that there always exists $k$ such that $J^{(k)}(R)$ is infinite and for this $k$, $I^{(k)}(r)\asymp J^{(k)}(r)$.
From this we deduce as above that the $k$th derivative of the Green function is infinite at $R$ and that its behavior is controlled by the $k$th derivatives of the Green functions associated with the first return kernels to parabolic subgroups.

In the case where parabolic subgroups are virtually abelian, it was proved in \cite{DPT} that the $k$th derivatives of the Green functions of first return kernels have a prescribed behavior.
This is because the first return kernel satisfies its own local limit theorem of the form $p_n\sim C\rho^nn^{-d/2}$, where $d$ is the rank of the parabolic subgroup.
Then, because of the above explanations, the $k$th derivative of the initial Green function associated with $\mu$ has the same behavior, 
and using Tauberian theory, one can deduce that $\mu^{(n)}(e)\sim C R^{-n}n^{-d/2}$.

Whenever parabolic subgroups are virtually nilpotent, several difficulties come up with this strategy.
The main problem is that the first return kernel to a parabolic subgroup is not Markov and does not have finite support.
It has finite exponential moments, but it is not enough to directly deduce that it satisfies its own local limit theorem.
Moreover, the first return kernel depends on $r$.
Our strategy is to refine a comparison theorem of Pittet and Saloff-Coste \cite{PSC} and to prove that the $n$th convolution power of the first return kernel behaves like the $n$th convolution power of a simple random walk on the parabolic subgroup, see Proposition~\ref{p:applicationsonctinuousPSCfirstreturn}.

In fact, although it is not part of Theorem~\ref{t:maintheoremsummary}, this refinement of Pittet-Saloff-Coste's theorem is one of our main results.
It can be roughly stated as follows: the behavior of the convolution powers of a continuous family of Markov transition kernels on a finite extension of an amenable finitely generated group is stable; see Proposition~\ref{p:continuousPSC} and Corollary~\ref{c:continuousPSCdouble} for precise statements.

\medskip
This concludes the proof in the convergent case.
Next, in the case where both $I^{(1)}(r)$ and $J^{(2)}(r)$ are infinite at $R$ (Theorem~\ref{t:maintheoremcritical}), we see that these two quantities play a role in the fundamental asymptotic differential equation~(\ref{e:fundamental}).
This is interesting and is reminiscent of the idea that both the underlying hyperbolic structure and the parabolic subgroups are influential.
Using again our continuous version of Pittet and Saloff-Coste comparison, we get that $J^{(2)}(r)$ has a prescribed behavior, as in the convergent case.
This in turn yields a new asymptotic differential equation satisfied by the Green function of $\mu$, which actually depends on the parity of the rank of spectral degeneracy, see~(\ref{e:newdifferentialequationcriticalcasedeven}) and~(\ref{e:newdifferentialequationcriticalcasedodd}).
From this we again derive the behavior of $\mu^{(n)}(e)$ using Tauberian theory.

As explained above, the asymptotics in Theorem~\ref{t:maintheoremcritical} were already obtained in one very specific example in \cite{DPT23}.
The proof here follows the same overall strategy, but some arguments were simpler in \cite{DPT23}, due to the combinatorial structure of free products.
The main intermediate result that had to be proved is given by Proposition~\ref{p:rho'G'}.
It asserts that the derivative of the spectral radius of the first return kernel to a parabolic subgroup behaves like the derivative of the Green function.
This is new and holds for every amenable parabolic subgroup.

\medskip
To conclude this introduction, let us briefly explain how the paper is organized.
In Section~\ref{s:background}, we review some material that will be needed throughout the paper. Namely, we recall the definition of relatively hyperbolic groups, we recall the local limit theorem for random walks on nilpotent groups and we state some results of the papers \cite{DussauleLLT1}, \cite{DussauleLLT2} and \cite{DGstability} that will be used afterward.

In Section~\ref{s:PSC}, we prove our generalization of Pittet-Saloff-Coste's Comparison Theorem for continuous families of Markov transition kernels on finite extensions of amenable groups.
Then in Section~\ref{s:applicationsPSC}, we apply this to the first return kernel to a parabolic subgroup.

Section~\ref{s:derivativespectralradius} is again an intermediate and technical section.
It is dedicated to proving that the derivative of the spectral radius of the first return kernel to a parabolic subgroup behaves like $I^{(1)}(r)$ as $r$ tends to $R$.
This will be important within the proofs of the local limit theorems.

The next two sections are the core of the paper, where we prove, respectively, Theorem~\ref{t:maintheoremconvergent} and Theorem~\ref{t:maintheoremcritical}.
We make the additional assumption that the random walk is lazy, that is, $\mu(e)>0$.
This is mainly because such an assumption is used in Section~\ref{s:applicationsPSC}, as it greatly simplifies the arguments.
We explain how to get rid of this assumption in the short Section~\ref{s:lazynotlazy}.
We also explain there how to treat the case where the random walk is not aperiodic.

Finally, in Section~\ref{s:conclusionandclassification}, we summarize the classification of local limit theorems on relatively hyperbolic groups, as it is known so far.
We also make further remarks and explain what is still missing to obtain a full classification beyond the case of virtually nilpotent parabolic subgroups.

\section{Some background and terminology}\label{s:background}
\subsection{Relatively hyperbolic groups}\label{ss:defrh}
Relatively hyperbolic groups have been considered by many authors with several equivalent points of view.
We briefly recall here their definition, using the terminology of Bowditch \cite{Bowditch}, Farb~\cite{Farb} and Osin \cite{Osin06}.

Consider a discrete group $\Gamma$ acting by isometries on a Gromov-hyperbolic space $X$. Let $o\in X$ be fixed. Define the limit set $\Lambda \Gamma$ as the closure of $\Gamma o$ in the Gromov boundary $\partial X$ of $X$. This set does not depend on $o$.

A point $\xi\in \Lambda \Gamma$ is called conical if there is a sequence $g_{n}$ in  $\Gamma$ and distinct points $\xi_1, \xi_2$ in $\Lambda \Gamma$ such that,  for  all $\xi\neq \zeta$ in $\Lambda \Gamma$,  the sequences 
$g_n\xi$  and $g_n\zeta$ converge to $\xi_1$ and  $\xi_2$ respectively.
A point $\xi\in \Lambda \Gamma$ is called parabolic if its stabilizer $\Gamma_\xi$  in $\Gamma$ is infinite and the elements of $\Gamma_\xi$  only fix $\xi$ in $\Lambda\Gamma$.
A parabolic limit point $\xi$ in $\Lambda \Gamma$ is said to be bounded parabolic if $\Gamma_\xi$ acts cocompactly on $\Lambda \Gamma \setminus \{\xi\}$. The action of $\Gamma$ on $X$ is called geometrically finite if $\Lambda \Gamma$  only contains conical limit points and bounded parabolic limit points.


Let $\Gamma$ be a finitely generated group and $S$ be a fixed generating set. 
Let $\Omega_0$ be a finite collection of subgroups,  none of them being conjugate. Let $\Omega$ be the closure of $\Omega_0$ under conjugacy.

The relative graph $\hat \Gamma = \hat\Gamma(S,  \Omega_0)$  is the Cayley graph of $\Gamma$ with respect to $S$ and the union of all $\mathcal{H}\in \Omega_0$.
It was introduced in \cite{Osin06} and
it is quasi-isometric to the coned-off graph introduced by Farb in \cite{Farb}.
The distance $\hat{d}$ in $\hat \Gamma$ is called the relative distance. We also denote by $\hat{S}_n$ the sphere of radius $n$ centered at $e$ in $\hat{\Gamma}$. Also, a relative geodesic is a geodesic in $\hat{\Gamma}$.

\begin{theorem}[\cite{Bowditch}]
Using the previous notations, the following conditions are equivalent.

\begin{enumerate}
\item The group $\Gamma$ has a geometrically finite action on a Gromov hyperbolic space $X$ such that the parabolic limit points are exactly the fixed points of elements in $\Omega$.

\item The relative graph $\hat \Gamma(S,  \Omega_0)$ is Gromov hyperbolic for the relative distance $\hat d$,  and for all $L>0$ and all $x\in \hat \Gamma$,  there exist finitely many closed loops of length $L>0$ which contain $x$.
\end{enumerate}
\end{theorem}

When these conditions are satisfied,  the group $\Gamma$ is said to be relatively hyperbolic with respect to $\Omega$.
We also say that $\Gamma$ is relatively hyperbolic with respect to $\Omega_0$, if $\Omega_0$ is fixed.
The elements of $\Omega$ are called (maximal) parabolic subgroups.

Assume now that $\Gamma$ is relatively hyperbolic with respect to $\Omega$,  and let $X$ be a Gromov hyperbolic space on which $\Gamma$ has a geometrically finite action whose parabolic subgroups  are the element of $\Omega$. The limit set $\Lambda \Gamma\subset \partial X$ is called the Bowditch boundary of $\Gamma$. It is unique up to equivariant homeomorphism.
The group $\Gamma$ is called non-elementary if its Bowditch boundary is infinite.

\medskip
Let $\Gamma=H_0*H_1$ be a free product of finitely generated groups.
Then, $\Gamma$ is relatively hyperbolic with respect to $H_0$ and $H_1$.
In fact, free products are combinatorial models for relatively hyperbolic groups.
In this example (assuming further that $H_0$ and $H_1$ are infinite), the Bowditch boundary is obtained by gluing together one point at infinity for each coset $gH_0$ and for each coset $gH_1$ with the set of infinite words, that is, infinite sequences $a_1b_1...a_nb_n...$ where $a_i\in H_0$, $b_i\in H_1$, $b_i\neq e$ and $a_i\neq e$ except maybe $a_1$.

Other archetypal examples are fundamental groups of finite-volume manifolds with pinched negative curvature.
Let $\mathcal{M}$ be such a manifold and let $\widetilde{\mathcal{M}}$ be its universal cover.
Then, the action of $\pi_1(\mathcal{M})$ on $\widetilde{\mathcal{M}}$ is geometrically finite.
To each cusp of $\mathcal{M}$ corresponds one orbit of parabolic limit points in the Gromov boundary of $\widetilde{\mathcal{M}}$.
The group $\pi_1(\mathcal{M})$ is relatively hyperbolic with respect to the stabilizers of the parabolic limit points for this action, which are also called the cusp subgroups.
The Bowditch boundary is the Gromov boundary of $\widetilde{\mathcal{M}}$.
See \cite{Bowditch95} for more details.
In fact, in this setting, the parabolic subgroups are necessarily virtually nilpotent by Margulis Lemma, which gives motivation for our classification.

\subsection{Virtually nilpotent groups}\label{ss:nilpotentgeometry}
Let $\Gamma$ be a group.
Given two subgroups $H_1$ and $H_2$, one defines $[H_1,H_2]$ as the subgroup of $\Gamma$ generated by all commutators of the form $xyx^{-1}y^{-1}$, $x\in H_1,y\in H_2$.
Then, one sets $\Gamma^{(1)}=\Gamma$, $\Gamma^{(2)}=[\Gamma,\Gamma]$ and defines inductively $\Gamma^{(n)}=[\Gamma^{(n-1)},\Gamma]$.
Since all subgroups $\Gamma^{(j)}$ are normal, $\Gamma^{(j)}$ is a subgroup of $\Gamma^{(j-1)}$.
A group is called nilpotent if there exists $n$ such that $\Gamma^{(n)}$ is the trivial subgroup.

Let $\Gamma$ be a finitely generated nilpotent group.
Denote by $N_\Gamma$ the largest integer such that $\Gamma^{(n)}$ is not trivial.
Then, $N_\Gamma$ is called the nilpotency step of $\Gamma$.
Note that all groups $\Gamma^{(j)}/\Gamma^{(j+1)}$ are finitely generated abelian groups, hence they have a well-defined rank $d_j$.
The homogeneous dimension of $\Gamma$ is then defined as
$$d_\Gamma=\sum_{j=1}^{N_\Gamma}jd_j.$$

A finitely generated nilpotent group $\Gamma$ has polynomial growth of degree $d_\Gamma$, that is, for any finite generating set,
$\big|\{x,d(e,x)\leq n\}\big|\asymp n^{d_\Gamma}$.
In particular, the homogeneous dimension is invariant under quasi-isometry.
Therefore, given a finitely generated virtually nilpotent group $\Gamma$, one can define its homogeneous dimension as the homogeneous dimension of any finite index nilpotent subgroup.

The homogeneous dimension is also a key parameter with regard to random walks on nilpotent groups.
It appears in the following local limit theorem that will be used at several places in the paper.

\begin{proposition}\label{p:LLTnilpotent}
Consider a finitely generated virtually nilpotent group $\Gamma$ and let $\mu$ be an admissible symmetric probability measure with finite second moments.
Let $d$ be the homogeneous dimension of $\Gamma$.
Then,
$$p_n(e,e)\asymp n^{-d/2}.$$
\end{proposition}

\begin{proof}
This result is classical for probability measures with finite support; see, for instance \cite[Theorem~15.8~(a)]{Woess-book}.
In fact, Alexopoulos \cite{Alexopoulos} proved that in such case $p_n\sim Cn^{-d/2}$.
The result for probability measures with finite second moments then follows from Pittet-Saloff-Coste's Comparison Theorem \cite{PSC}, see also \cite[Theorem~15.1, Theorem~15.8]{Woess-book}.
\end{proof}

\subsection{The first return kernel to a parabolic subgroup}\label{ss:firstreturnkerneldef}
Let $\Gamma$ be a relatively hyperbolic group and let $\mu$ be a probability measure on $\Gamma$.
Given a subset $A$ of $\Gamma$, we may consider the Green function restricted to paths staying in $A$, except maybe at the endpoints, which is defined for any $x,y\in \Gamma$ by
$$G(x,y;A|r)=\sum_{n\geq 1}\sum_{\underset{z_1,...,z_{n-1}\in A}{z_0=x,z_1,...,z_n=y}}r^n\mu(z_0^{-1}z_1)\mu(z_1^{-1}z_2)...\mu(z_{n-1}^{-1}z_n).$$
Let $H$ be a parabolic subgroup. The first return kernel to $H$ associated with $r\mu$ is defined as
$$p_{r,H}(x,y)=G(x,y;H^c|r).$$
For $r=1$, $p_{r,H}(x,y)$ is the probability that the random walk, starting at $x$, eventually returns to $H$ and that the first return to $H$ is at $y$.

More generally, given $\eta\geq 0$, define similarly the first return kernel to the $\eta$-neighborhood $N_\eta(H)$ of $H$ as 
$$p_{r,H,\eta}(x,y)=G(x,y;N_\eta(H)^c|r).$$
Next, define the $n$th convolution power of $p_{r,H,\eta}$ as
$$p_{r,H,\eta}^{(n)}(x,y)=\sum_{z_1,...,z_{n-1}\in H}p_{r,H,\eta}(x,z_1)p_{r,H,\eta}(z_1,z_2)...p_{r,H,\eta}(z_{n-1},y)$$
and consider the associated Green function
$$G_{r,H,\eta}(x,y|t)=\sum_{n\geq 0}t^np_{r,H,\eta}^{(n)}(x,y).$$
Also, let $\rho_{H,\eta}(r)$ be the spectral radius of $p_{r,H,\eta}$ defined as
$$\rho_{H,\eta}(r)=\limsup \bigg(p_{r,H,\eta}^{(n)}(e,e)\bigg)^{1/n}.$$
Let $R_{H,\eta}(r)=\rho_{H,\eta}(r)^{-1}$ be its inverse.
Then, $R_{H,\eta}(r)$ is the radius of convergence of the Green function $G_{r,H,\eta}$.

\begin{lemma}\label{l:sameGreen}
Let $H$ be a parabolic subgroup, let $\eta\geq0$ and let $r\leq R$.
Then, for every $x,y\in N_\eta(H)$, we have
$$G_{r,H,\eta}(x,y|1)=G(x,y|r).$$
\end{lemma}

Indeed, the Green function $G(x,y|r)$ counts all trajectories from $x$ to $y$.
Then, $p_{r,H,\eta}^{(n)}(x,y)$ counts all trajectories with exactly $n$ steps in $H$.
By conditioning on successive visits to $H$, we get the desired identification.
We refer to \cite[Lemma~4.4]{DGstability} for a complete proof.

Since $\Gamma$ is non-amenable, $G(e,e|R)$ is finite by \cite[Theorem~7.8]{Woess-book}.
This implies in particular that $G_{R,H,\eta}(e,e|1)$ is finite.
Therefore, $R_{H,\eta}(R)\geq 1$.
The following notion was introduced in \cite{DGstability}.

\begin{definition}\label{d:spectraldegeneracy}
Say that $\mu$ is spectrally degenerate along a parabolic subgroup $H$ if $\rho_{H,0}(R)=1$.
\end{definition}

Note that by \cite[Lemma~4.9]{DGstability}, if $\mu$ is spectrally degenerate along $H$, then for every $\eta\geq 0$, $\rho_{H,\eta}(R)=1$.
In other words, if $\mu$ is spectrally degenerate along $H$, it is also spectrally degenerate along any neighborhood of $H$.
This notion was first introduced to study the homeomorphism type of the $R$-Martin boundary in \cite{DGstability}, but it turned out to be also important in the study of local limit theorems \cite{DussauleLLT2}.

\medskip
A large part of this paper is dedicated to studying the asymptotics properties of $p_{r,H,\eta}^{(n)}$.
For this, we will often use the following identification.
Note that $H$ acts on $N_\eta(H)$ by left multiplication.
Any section of $N_\eta(H)\to N_\eta(H)/H$ provides an identification of
$N_\eta(H)$ as $H\times E$, where $E$ is a finite set.
As a consequence, $p_{r,H,\eta}$ can be viewed as a transition kernel on $H\times E$ where $E$ is a finite set.
Since the initial random walk is invariant by left multiplication, $p_{r,H,\eta}$ is $H$-invariant, that is for any $h\in H$ and any $x,y\in N_\eta(H)$,
$p_{r,H,\eta}(x,y)=p_{r,H,\eta}(hx,hy)$.
However, note that $p_{r,H,\eta}$ is not Markov, that is, given $x\in N_\eta(H)$, the total mass of $p_{r,H,\eta}$ at $x$ is not constant equal to 1.

Finally, let us introduce the following notation.
An element of $N_\eta(H)$, identified with $H\times E$, can be written as $(x,i)$, $x\in H$, $i\in E$.
In the sequel, given a transition kernel on $H\times E$ we will often write
$$p_{i,j}(x,y)=p((x,i),(y,j)).$$
This will allow us to consider $p$ as a family of transition kernels on $H$ and to consider the matrix indexed by $E$, whose $(i,j)$ entry is given by $p_{i,j}(x,y)$.

\subsection{Derivatives of Green functions}\label{ss:derivativesGreenfunction}
Let $\Gamma$ be a finitely generated group.
Recall that the quantities $I^{(1)}(r)$ and $I^{(2)}(r)$ were defined in the Introduction.
More generally, we define $I^{(k)}(r)$ as follows:
\begin{equation}\label{e:defIk}
I^{(k)}(r)=\sum_{x_1,...,x_k\in \Gamma}G(e,x_1|r)G(x_1,x_2|r)...G(x_{k-1},x_k|r)G(x_k,e|r).
\end{equation}
Also, for $x\in \Gamma$, we set
\begin{equation}\label{e:defIkg}
I^{(k)}_x(r)=\sum_{x_1,...,x_k\in \Gamma}G(e,x_1|r)G(x_1,x_2|r)...G(x_{k-1},x_k|r)G(x_k,x|r).
\end{equation}

\begin{lemma}\label{l:GreenderivativeI1}\cite[Proposition~1.9]{GouezelLalley}
For every $x\in \Gamma$ for every $r\leq R$, we have
$$I^{(1)}_x(r)=\frac{d}{dr}\big (rG(e,x|r)\big)=rG'(e,x|r)+G(e,x|r).$$
\end{lemma}

Similar expressions hold for higher derivatives; see \cite[Lemma~3.2]{DussauleLLT1}.
In particular, for the second derivative, we have
\begin{equation}\label{e:I2}
2I^{(2)}_x(r)=r^2G''(e,x|r)+4rG'(e,x|r)+2G(e,x|r).
\end{equation}

We now assume that $\Gamma$ is non-elementary relatively hyperbolic.
As explained above, since $\Gamma$ is non-amenable, $G(e,e|R)$ is finite.
In particular, we have

\begin{equation}\label{e:I1G1}
I^{(1)}(r)\asymp G'(e,e|r).
\end{equation}

Also, because of the fundamental asymptotic differential equation~(\ref{e:fundamental}),
we see that either $G''(e,e|R)$ is finite or that $I^{(2)}(r)/G'(e,e|r)$ tends to infinity, so that
\begin{equation}\label{e:I2G2}
I^{(2)}(r)\asymp G''(e,e|r).
\end{equation}

Given a parabolic subgroup $H$, $\eta\geq 0$ and $x\in N_\eta(H)$, we also introduce the quantity
\begin{equation}\label{e:defIHk}
I_{H,\eta,x}^{k}(r)=\sum_{x_1,...,x_k\in N_\eta(H)}G(e,x_1|r)G(x_1,x_2|r)...G(x_{k-1},x_k|r)G(x_k,x|r).
\end{equation}
We write $I^{(k)}_{H,\eta}(r)=I^{(k)}_{H,\eta,e}(r)$ for simplicity.
Choosing a finite set $\Omega_0$ of representatives of conjugacy classes of parabolic subgroups, we then set
\begin{equation}\label{e:defJeta}
J^{(k)}_\eta(r)=\sum_{H\in \Omega_0}I^{(k)}_{H,\eta}(r).
\end{equation}
For $\eta=0$, we simply write $J^{(k)}(r)=J_0^{(k)}(r)$.

Recall that $p_{r,H,\eta}$ is the transition kernel of first return to $N_\eta(H)$ and that $G_{r,H,\eta}$ is the associated Green function.
Also, recall that by Lemma~\ref{l:sameGreen}, for every $x,y\in N_\eta(H)$, $G(x,y|r)=G_{r,H,\eta}(x,y|1)$.
Consequently, we have
$$I_{H,\eta,x}^{(k)}(r)=\sum_{x_1,...,x_k\in N_\eta(H)}G_{r,H,\eta}(e,x_1|1)G_{r,H,\eta}(x_1,x_2|1)...G_{r,H,\eta}(x_k,x|1).$$
Therefore, $I_{H,\eta,x}^{(1)}(r)$ can be interpreted as the derivative of $t\mapsto tG_{r,H,\eta}(e,x|t)$ at $t=1$.
Once again, similar formulae hold for higher derivatives.
We deduce the following.
\begin{lemma}
Assuming that $\mu$ is not spectrally degenerate along $H$, for any $k$, $I_{H,\eta}^{(k)}(R)$ is finite.
\end{lemma}

In fact, we have the following important result.

\begin{proposition}\label{p:NSDimpliesSPR}\cite[Proposition~3.7]{DussauleLLT1}
Let $\Gamma$ be a relatively hyperbolic group and let $\mu$ be a finitely supported admissible symmetric probability measure on $\Gamma$.
If $\mu$ is not spectrally degenerate, then it is spectrally positive recurrent.
\end{proposition}

In view of what precedes, the nontrivial part of this proposition is that spectral non-degeneracy implies that $\mu$ is divergent.
We conclude this section with the following result.

\begin{proposition}\label{p:J1finite}\cite[Proposition~6.3]{DGstability}
Let $\Gamma$ be a relatively hyperbolic group and let $\mu$ be a finitely supported admissible symmetric probability measure on $\Gamma$.
For any parabolic subgroup $H$, for any $\eta\geq 0$ and for any $r\leq R$,
$I_{H,\eta}^{(1)}(r)$ is finite.
\end{proposition}

\section{A comparison theorem}\label{s:PSC}
In the following sections, we will need to deal with transition kernels on a fixed neighborhood of a parabolic subgroup $H$.
As seen above, such a neighborhood can be identified with a finite extension of $H$.
Here we work in the general framework of finite extension of finitely generated amenable groups.
We set $X=\Gamma\times E$ where $E$ is a finite set and where $\Gamma$ is a finitely generated amenable group.
An element of $X$ will be written $(x,j)$, with $x\in \Gamma$ and $j\in E$.
Given a transition kernel $p$ on $X$, we write $p_{i,j}(x,y)=p((x,i),(y,j))$.
We say that $p$ is a Markov transition kernel if for every $(x,i)\in \Gamma\times E$,
$$\sum_{(y,j)\in X}p_{i,j}(x,y)=1.$$

Our goal is to prove a continuous version of Pittet-Saloff-Coste's Comparison Theorem proved in \cite{PSC}.
Namely, the result of Pittet and Saloff-Coste states that the asymptotic type of the convolution powers of a symmetric Markov transition kernel on a finitely generated group is stable by quasi-isometry.
We generalize this result to a continuous family of strongly reversible Markov transition kernels; see Corollary~\ref{c:continuousPSCdouble}.
We will then apply this to the transition kernels of first return to parabolic subgroups in the next section.

\subsection{A trace theorem}
We consider a Markov transition kernel $p$ on $\Gamma\times E$.
Given a subset $A$ of $\Gamma \times E$, we write $q_A$ for the sub-Markov transition kernel restricted to $A$, defined by
$$
q_{A;i,j}(x,y):=q_A((x,i),(y,j))=\left \{
\begin{array}{cc}
 p_{i,j}(x,y)    &  \text{ if } (x,i),(y,j)\in A \\
 0    &  \text{ otherwise.}
\end{array}\right.$$
Whenever $A$ is finite, we let $Q_A$ be the matrix indexed by elements of $A$ whose $(x,i),(y,j)$ entry is given by $q_{A;i,j}(x,y)$.

Also, if $A$ is a subset of $\Gamma$, we set
$$\tilde A=\{(x,i)\in \Gamma\times E, x\in A\}=A\times E.$$
We prove the following.
\begin{proposition}\label{p:tracetheorem}
Assume that $\Gamma$ is an amenable finitely generated group and $E$ is a finite set.
Let $p$ be a $\Gamma$-invariant Markov transition kernel on $\Gamma\times E$.
Then
$$\frac{1}{|E|}\sum_{i\in E}p_{i,i}^{(n)}(e,e)=\sup \left \{\frac{1}{|\tilde A|}\mathrm{Tr}(Q_{\tilde A}^n),A\text{ finite subset of }\Gamma\right \}.$$
Moreover,
$$\frac{1}{|E|}\sum_{i\in E}p_{i,i}^{(n)}(e,e)\leq \sup \left \{\frac{1}{|B|}\mathrm{Tr}(Q_{B}^n),B\text{ finite subset of }\Gamma\times E\right \}\leq \sum_{i\in E}p_{i,i}^{(n)}(e,e).$$
\end{proposition}

\begin{proof}
We consider a finite set $K_\epsilon\subset \Gamma$ containing $e$ such that for every $i$,
$$\sum_j\sum_{x\in K_\epsilon}p_{i,j}(e,x)\geq 1-\epsilon.$$
Second, fixing $n$, by F{\o}lner criterion, we consider a sequence of finite sets $U_k\subset \Gamma$ such that
$$\frac{|U_kK_\epsilon^n|}{|U_k|}\underset{k\to \infty}{\longrightarrow}1.$$
We set $A_k=U_kK_\epsilon^n$.
We write $X_k=(x_k,i_k)$ the Markov chain at time $k$.
Then, for $x\in U_k$, letting $q$ be the sub-Markov transition kernel restricted to $\tilde A_k$,
$$q_{i,i}^{(n)}(x,x)=\mathbb P_{x,i}(X_n=(x,i),X_j\in \tilde A_k,j=1,...,n-1).$$
Thus,
$$q_{i,i}^{(n)}(x,x)=\mathbb P_{x,i}(X_n=(x,i))-\mathbb P_{x,i}(X_n=(x,i),\exists j, X_j\notin \tilde A_k).$$
Next, if for all $1\leq j\leq n$, $x_{j-1}^{-1}x_j\in K_\epsilon$, where $x_0=x$, then for all $1\leq j\leq n$,
$x_j=x(x_0^{-1}x_1)(x_1^{-1}x_2)...(x_{j-1}^{-1}x_j)\in xK_\epsilon^j\subset A_k$, since $x\in U_k$ and $e\in K_\epsilon$, hence $X_j\in \tilde A_k$.
Consequently,
$$\mathbb P_{x,i}(X_n=(x,i),\exists j, X_j\notin \tilde A_k)\leq \mathbb P_{x,i}(\exists j, x_{j-1}^{-1}x_j\notin K_\epsilon)\leq \sum_{j=1}^n\mathbb P_{x,i}(x_j\notin x_{j-1}K_\epsilon).$$
Given $1\leq j\leq n$,
$$\mathbb P_{x,i}(x_j\notin x_{j-1}K_\epsilon)=\sum_{y,l}\mathbb P_{x,i}(X_{j-1}=(y,l))\mathbb P_{x,i}(x_j\notin yK_\epsilon|X_{j-1}=(y,l)).$$
Finally, by the Markov property,
$$\mathbb P_{x,i}(x_j\notin x_{j-1}K_\epsilon)=\sum_{y,l}\sum_m\sum_{z\notin yK_\epsilon}\mathbb P_{x,i}(X_{j-1}=(y,l))p_{l,m}(y,z)\leq \epsilon.$$
We thus get
$$\mathbb P_{x,i}(X_n=(x,i),\exists j, X_j\notin \tilde A_k)\leq n\epsilon.$$
We conclude that
$$q_{i,i}^{(n)}(x,x)\geq p_{i,i}^{(n)}(e,e)-n\epsilon.$$
Summing over $i\in E$ and $x\in U_k$,
$$\sum_{i}p_{i,i}^{(n)}(e,e)-|E|n\epsilon\leq \frac{1}{|U_k|}\mathrm{Tr}(Q_{\tilde A_k}^n)\leq \frac{|\tilde A_k|}{|U_k|}\frac{1}{|\tilde A_k|}\mathrm{Tr}(Q_{\tilde A_k}^n).$$
Note that $|\tilde A_k|=|E||A_k|$, so that $|\tilde A_k|/|U_k|$ converges to $|E|$ as $k$ tends to infinity.
Letting $k$ tend to infinity and $\epsilon$ tend to 0, we get
$$\frac{1}{|E|}\sum_{i}p_{i,i}^{(n)}(e,e)\leq \sup \left \{\frac{1}{|\tilde A|}\mathrm{Tr}(Q_{\tilde A}^n),A \text{ finite subset of } \Gamma\right \}.$$
This proves the first inequality.
We also directly get that
$$\frac{1}{|E|}\sum_{i}p_{i,i}^{(n)}(e,e)\leq \sup \left \{\frac{1}{| B|}\mathrm{Tr}(Q_{B}^n), B\text{ finite subset of } \Gamma\times E\right \}.$$

Next, consider a finite subset $B$ of $\Gamma\times E$.
For every $i\in E$, we set
$$B(i)=\{x\in \Gamma,(x,i)\in B\}$$ which is a subset of $\Gamma$.
For every $i,j$ and every $x,y$,
$$q^{(n)}_{B;i,j}(x,y)\leq p^{(n)}_{i,j}(x,y).$$
Therefore,
$$\mathrm{Tr}(Q_B^n)\leq \sum_{x,i}p^{(n)}_{i,i}(x,x)=\sum_{i}|B(i)|p^{(n)}_{i,i}(e,e).$$
First of all, note that $|B(i)|\leq |B|$, so
$$\mathrm{Tr}(Q_B^n)\leq |B|\sum_{i}p^{(n)}_{i,i}(e,e).$$
Second, if $A$ is a finite subset of $\Gamma$ and $B=\tilde A$, then for every $i\in E$, we have $|B(i)|=|A|=|\tilde A|/|E|$.
Thus,
$$\mathrm{Tr}(Q_{\tilde A}^n)\leq |\tilde A|\frac{1}{|E|}\sum_{i}p^{(n)}_{i,i}(e,e).$$
This proves the two last inequalities.
\end{proof}

\subsection{Continuous families of Markov transition kernels}
In preparation of our comparison theorem, we introduce the following set of definitions.
\begin{definition}
A continuous family of Markov transition kernels on $X$ is a map
$$r\in I\mapsto p(r),$$
where $I$ is an interval of $\mathbb R$, $p(r)$ is a Markov transition kernel on $X$ for every $r$ and for every $(x,i),(y,j)$ in $X$, the map
$$r\mapsto p(r)((x,i),(y,j))$$
is continuous.
\end{definition}

For simplicity, we write $p_{r;i,j}(x,y)=p(r)((x,i),(y,j))$ in what follows.

\begin{definition}
A continuous family of Markov transition kernels on $X$ is said to be $\Gamma$-invariant if for every $r\in I$, for every $(x,i),(y,j)$ in $X$ and for every $g\in \Gamma$
$$p_{r;i,j}(gx,gy)=p_{r;i,j}(x,y).$$
\end{definition}

\begin{definition}
A continuous family of Markov transition kernels on $X$ is said to be symmetric if for every $r\in I$, for every $(x,i),(y,j)$,
$$p_{r;i,j}(x,y)=p_{r;j,i}(y,x).$$
\end{definition}

\begin{remark}
Note that $p(r)$ is asked to be symmetric, but given $i,j\in E$, $p_{r;i,j}$ might not be symmetric, except if $p_{r;i,j}=p_{r;j,i}$.
\end{remark}

\begin{definition}
A continuous family of Markov transition kernels $r\mapsto p(r)$ on $X$ is said to be reversible if there exists a map $r\mapsto m(r)$
where $m(r):X\to \mathbb R_+$ is a function satisfying the following properties.
Setting $m_{r;i}(x)=m(r)((x,i))$,
\begin{itemize}
    \item for every $(x,i)$ in $X$, the function $r\mapsto m_{r;i}(x)$ is continuous
    \item for every $r\in I$, for every $(x,i),(y,j)$ in $X$,
    $$m_{r;i}(x)p_{r;i,j}(x,y)=p_{r;j,i}(y,x)m_{r;j}(y).$$
\end{itemize}
The function $m$ is called the total conductance function.
The continuous family of Markov transition kernels $r\mapsto p(r)$ is called strongly reversible if $m$ is uniformly bounded from below and from above by continuous functions, i.e.\ there exist continuous positive functions $r\mapsto c(r)$ and $r\mapsto C(r)$ such that for every $r\in I$, for every $(x,i)$ in $X$,
$$c(r)\leq m_{r;i}(x) \leq C(r).$$
\end{definition}

We fix a finite generating set $S$ on $\Gamma$ and we write $x\sim y$ if $x^{-1}y\in S$.
Also, we write $p^{(K)}(r)$ for the $K$th power of convolution of $p(r)$ and for simplicity, we write
$p^{(K)}_{r;i,j}(x,y)=p^{(K)}(r)((x,i),(y,j))$.

\begin{definition}
A continuous family of Markov transition kernels on $X$ is said to be uniformly irreducible (with respect to the finite generating set $S$) if there exist a continuous function $r\mapsto \epsilon(r)$ and an integer $K$ such that for every $i,j$ in $E$ and for every $x,y$ such that $x\sim y$,
$$p^{(K)}_{r;i,j}(x,y)\geq \epsilon(r).$$
\end{definition}

\begin{remark}
If $p$ is a continuous family of Markov transition kernels which is $\Gamma$-invariant, then it is uniformly irreducible if and only if for every $x\in S$, for every $i,j$ in $E$,
$p_{r;i,j}^{(K)}(e,x)\geq \epsilon(r)$.
\end{remark}

We let $d$ be the word distance associated with $S$.
We say that a transition kernel $p$ has finite second moment if for every $i$ and for every $x$,
$$\sum_{y\in \Gamma}\sum_{j\in E}d(x,y)^2p_{i,j}(x,y)$$
is finite.
Whenever $p$ is $\Gamma$-invariant, this is equivalent to
$$\sum_{x\in \Gamma}\sum_{j\in E}d(e,x)^2p_{i,j}(e,x)$$ for every $i$.
We then set
$$\mathcal{M}_2(p)=\sum_{i\in E}\sum_{x\in \Gamma}\sum_{j\in E}d(e,x)^2p_{i,j}(e,x).$$

\begin{definition}
A continuous family of Markov transition kernels on $X$ is said to have continuous second moment (with respect to the finite generating set $S$) if for every $r\in I$, $p(r)$ has finite second moment $\mathcal{M}_2(r):=\mathcal{M}_2(p(r))$ and the function
$r\mapsto \mathcal{M}_2(r)$ is continuous.
\end{definition}

\subsection{Dirichlet forms}
Let $X$ be a graph and let $p$ be a transition kernel, which is reversible with respect to a total conductance function $m$.
The space $\ell_2(X,m)$ is the set of functions $f:X\to \mathbb R$ such that
$$\sum_{x\in X}f(x)^2m(x)<+\infty.$$
It can be endowed with the inner product
$$\langle f,g\rangle=\sum_{x\in X}f(x)g(x)m(x).$$
The Dirichlet form associated with $p$ is defined as
$$\mathcal{D}_p(f)=\frac{1}{2}\sum_{x,y\in X}\big(f(x)-f(y)\big)^2m(x)p(x,y).$$
Note that this quantity is called Dirichlet norm at several places in \cite{Woess-book}.
In fact, $\sqrt{\mathcal{D}_p(\cdot)}$ is a semi-norm whose kernel consists of constant functions.
In any case, $\mathcal{D}_p(f)$ is defined for any $f\in \ell_2(X,m)$.
In fact, one has
$$\mathcal{D}_p(f)=\langle \nabla f,\nabla f\rangle,$$
where $\nabla$ is the so-called difference operator.
We refer to \cite[Chapter~2.A]{Woess-book} for more details.

Assuming that $p$ is a Markov transition kernel, one has
\begin{equation}\label{e:fundamentalequationDirichletnorms}
\mathcal{D}_p(f)=\langle f,(I-P)f\rangle,
\end{equation}
where $P$ is the Markov operator defined by
$$P(f)(x)=\sum_{y\in X}p(x,y)f(y).$$
See for instance \cite[(2.2)]{Woess-book} for a proof.

Now consider, as above, a graph $X=\Gamma\times E$ where $\Gamma$ is a finitely generated amenable group and $E$ is finite.
Let $p$ be a reversible transition kernel on $X$.
The Dirichlet form can be written as
\begin{equation}\label{e:defDirichletnorm}
\mathcal{D}_p(f)=\frac{1}{2}\sum_{(x,i),(y,j)\in X}\big(f((x,i))-f((y,j))\big)^2m_i(x)p_{i,j}(x,y).
\end{equation}

Our next goal is to compare the Dirichlet form of a continuous family of Markov transition kernels with the Dirichlet form of a simple random walk on $\Gamma\times E$.

Fix a finite generating set $S$ on $\Gamma$ and let $\mu_S$ be the uniform measure on $S$.
Define the transition kernel $p_S$ on $X= \Gamma\times E$ by
$$p_S((x,i)(y,j))=\frac{1}{|E|}\mu_S(x^{-1}y).$$

In the remainder of this section, we consider a continuous family of Markov transition kernels $r\in I\mapsto p(r)$ which is reversible with respect to some total conductance function $m$.
For simplicity, we write
$\mathcal{D}_r=\mathcal{D}_{p(r)}$
and
$\mathcal{D}=\mathcal{D}_{p_S}$.
Note that $p_S$ is symmetric and that
\begin{equation}\label{e:DirichletnormSRW}
\mathcal{D}(f)=\frac{1}{2}\sum_{(x,i),(y,j)\in X}\big(f((x,i))-f((y,j))\big)^2\frac{1}{|E|}\frac{1}{|S|}.
\end{equation}

We denote by $\ell_0(X)$ the set of functions with finite supports on $X=\Gamma\times E$.

\begin{lemma}\label{l:epsilon0Dirichletnorms}
Assume that $r\mapsto p(r)$ is uniformly irreducible and strongly reversible.
Then, there exists a positive continuous function $\epsilon_0(r)$ such that for every $r\in I$ and for every $f\in \ell_0(X)$,
$$\mathcal{D}_{r}(f)\geq \epsilon_0(r) \mathcal{D}(f).$$
\end{lemma}

\begin{proof}
Since $r\mapsto p(r)$ is uniformly continuously irreducible, there exist $K$ and a function $\epsilon(r)$ such that
for $x\sim y$ and for every $i,j$, one has
$$p^{(K)}_{r;i,j}(x,y)\geq \epsilon(r)\geq \epsilon(r)\frac{1}{|E||S|}.$$
Moreover, the total conductance function $m$ satisfies $m_{r;i}(x)\geq c(r)$ for some positive continuous function $c$.
Consequently, by~(\ref{e:defDirichletnorm}) and~(\ref{e:DirichletnormSRW}),
$$\mathcal D_{p^{(K)}(r)}(f)\geq \epsilon(r)c(r)\mathcal D(f).$$
Now by \cite[Lemma~2.5]{Woess-book}, $\mathcal D_{p^{(K)}(r)}(f)\leq K^2\mathcal D_{r}(f)$, hence setting $\epsilon_0(r)=\frac{\epsilon(r)c(r)}{K^2}$,
the function $\epsilon_0$ is continuous, positive and satisfies $\mathcal{D}_{r}(f)\geq \epsilon_0(r) \mathcal{D}(f)$.
\end{proof}

\begin{lemma}\label{l:epsilon1Dirichletnorms}
Assume that $r\mapsto p(r)$ is $\Gamma$-invariant, strongly reversible and has continuous second moment.
Then, there exists a positive continuous function $\epsilon_1(r)$ such that for every $r\in I$, for every $f\in \ell_0(X)$,
$$\mathcal{D}(f)\geq \epsilon_1(r) \mathcal{D}_{r}(f).$$
\end{lemma}

\begin{proof}
Consider the function $\Phi$ defined on the set of edges of the graph $X$ by
$$\Phi(\mathfrak e)=\frac{1}{2}\sum_{(x,i),(y,i)}m_{r;i}(x)p_{r;i,j}(x,y)d(x,y)\frac{|\Pi_{\mathfrak e}((x,i),(y,j))|}{|\Pi((x,i),(y,j))|},$$
where $|\Pi((x,i),(y,j))|$ is the number of geodesics from $(x,i)$ to $(y,j)$ in the graph $X$ and $|\Pi_{\mathfrak e}((x,i),(y,j))|$ is the number of such geodesics containing the edge $\mathfrak e$.
It is proved in \cite[Theorem~3.2]{Woess-book} that this function $\Phi$ satisfies
$$\mathcal{D}_{r}(f)\leq \sup_{\mathfrak e} \Phi(\mathfrak e)\mathcal{D}(f).$$
Assuming further that $r\mapsto p(r)$ is $\Gamma$-invariant and that $m_{r;i}(x)\leq C(r)$ for some continuous function $C(r)$, it is proved in \cite[Proposition~3.20]{Woess-book} that the function $\Phi$ satisfies that for any $\mathfrak e$,
$$\Phi(\mathfrak e)\leq C(r)\mathcal{M}_2(r).$$
Therefore, setting $\epsilon_1(r)=\big(C(r)\mathcal{M}_2(r)\big)^{-1}$, the function $\epsilon_1$ is continuous, positive and satisfies
$\mathcal{D}(f)\geq \epsilon_1(r) \mathcal{D}_{r}(f)$.
\end{proof}

\subsection{A continuous version of Pittet-Saloff-Coste's Comparison Theorem}
Gathering all results of the previous subsections, we are ready to prove a continuous version of Pittet-Saloff-Coste's theorem for Markov chains on extensions of amenable groups.
\begin{proposition}\label{p:continuousPSC}
Let $\Gamma$ be a finitely generated amenable group and let $E$ be a finite set.
Let $r\mapsto p(r)$ be a continuous family of Markov transition kernels on $X=\Gamma\times E$.
Assume that it is $\Gamma$-invariant, strongly reversible and uniformly irreducible.
Consider a finite generating set $S$ for $\Gamma$ and let $\mu_S$ be the uniform measure on $S$.
There is a continuous function $r\mapsto \epsilon(r)\in (0,1/2)$ such that for every $n$,
$$\frac{1}{|E|}\sum_{i\in E}p^{(2n+2)}_{r;i,i}(e,e)\leq 2(1-\epsilon(r))^{2n}+\frac{2}{|E|}\mu_S^{(2\lfloor \epsilon(r)n\rfloor)}(e).$$
Moreover, assuming that $r\mapsto p(r)$ has continuous second moment, there exists a continuous function $r\mapsto \epsilon'(r)\in (0,1/2)$ such that
$$\frac{1}{|E|}\mu_S^{(2n+2)}(e)\leq 2(1-\epsilon'(r))^{2n}+\frac{2}{|E|}\sum_{i\in E}p^{(2\lfloor \epsilon'(r)n\rfloor)}_{r;i,i}(e,e).$$
\end{proposition}

\begin{proof}
We write $\langle \cdot,\cdot\rangle_m$ for the scalar product
$$\langle f,g\rangle_m=\sum_{(x,i)\in X}f((x,i))g((x,i))m_{r;i}(x)$$
and we write $\langle\cdot,\cdot\rangle$ for the standard scalar product
$$\langle f,g\rangle=\sum_{(x,i)\in X}f((x,i))g((x,i)).$$
Since $r\mapsto p(r)$ is strongly reversible, we have for any $f\in \ell_0(X)$
$$\langle f,f\rangle_m\leq C(r) \langle f,f\rangle.$$
We deduce from Lemma~\ref{l:epsilon0Dirichletnorms} that
$$\frac{\mathcal{D}_r(f)}{\langle f,f\rangle_m}\geq \frac{\epsilon_0(r)}{C(r)}\frac{\mathcal{D}(f)}{\langle f,f\rangle}.$$
We set
$$\epsilon(r)=\min\left\{\frac{1}{2},\frac{\epsilon_0(r)}{2C(r)}\right \}.$$
Then, $r\mapsto \epsilon(r)$ is continuous, $\epsilon(r)\in (0,1/2)$ and it satisfies
\begin{equation}\label{e:comparisonminmaxepsilon}
\frac{\mathcal{D}_r(f)}{\langle f,f\rangle_m}\geq 2\epsilon(r)\frac{\mathcal{D}(f)}{\langle f,f\rangle}.
\end{equation}
We now fix a finite set $A$ of $\Gamma$ and consider the subset $\tilde A=A\times E$ of $X$.
We let $Q_r$ be the matrix indexed by elements of $\tilde A$, whose $((x,i),(y,j))$ entry is given by
$$
Q_r((x,i),(y,j))=\left \{
\begin{array}{cc}
 p_{r;i,j}(x,y)   &  \text{ if } x,y\in A \\
 0    &  \text{ otherwise.}
\end{array}\right.$$
We also let $Q$ be the matrix whose $((x,i),(y,j))$ entry is given by
$$
Q((x,i),(y,j))=\left \{
\begin{array}{cc}
\frac{1}{|E|}\mu_S(x^{-1}y)   &  \text{ if } x,y\in A \\
 0    &  \text{ otherwise.}
\end{array}\right.$$
Then, for any $f$ whose support lies in $\tilde A$, we have
$$\mathcal{D}_r(f)=\langle f,(I-Q_r)f\rangle_m$$
and
$$\mathcal{D}(f)=\langle f,(I-Q)f\rangle.$$
The matrix $Q_r$, respectively the matrix $Q$, is self-adjoint with respect to the scalar product $\langle \cdot,\cdot\rangle_m$, respectively the scalar product $\langle\cdot,\cdot\rangle$.
We let $\lambda_r(1)\geq ... \geq \lambda_r(|\tilde A|)$ be the eigenvalues of $Q_r$ and
$\lambda(1)\geq ... \geq \lambda(|\tilde A|)$ be those of $Q$.
The min-max theorem shows that
$$1-\lambda_r(k)=\min_{\mathrm{dim}V=k}\max_{f\in V}\frac{\langle f,(I-Q_r)f\rangle_m}{\langle f,f\rangle_m}$$
and similarly,
$$1-\lambda(k)=\min_{\mathrm{dim}V=k}\max_{f\in V}\frac{\langle f,(I-Q)f\rangle}{\langle f,f\rangle}$$
where the minimum is over subspaces $V$ of $\ell_2(\tilde A)$.
We thus deduce from~(\ref{e:comparisonminmaxepsilon}) that for every $1\leq k\leq |\tilde A|$,
$$1-\lambda_r(k)\geq 2\epsilon(r)(1-\lambda(k)).$$
See also \cite[(2.3)]{DSC} for a similar statement.
Therefore, if $\lambda(k)\geq 1/2$, using that $1-u\leq \mathrm{e}^{-u}$ for $u\geq 0$ and that $\mathrm{e}^{v-1}\leq \sqrt{v}$ for $1/2\leq v\leq 1$, we get
$$\lambda_r(k)\leq 1-2\epsilon(r)(1-\lambda(k))\leq \mathrm{e}^{-2\epsilon(r)(1-\lambda(k))}\leq \lambda(k)^{\epsilon(r)}.$$
Assuming further that $\lambda_r(k)\geq 0$, we get
\begin{equation}\label{e:comparingeigenvalueslambdakgeq1/2}
\lambda_r(k)^{2n}\leq \lambda(k)^{2n\epsilon(r)}.
\end{equation}
On the other hand, if $\lambda(k)<1/2$, then
\begin{equation}\label{e:comparingeigenvalueslambdakleq1/2}
\lambda_r(k)\leq (1-\epsilon(r)).
\end{equation}
Next, observe that the trace of $Q_r^{2n+1}$ is non-negative, as
$$\mathrm{Tr}(Q_r^{2n+1})=\sum_{(x,i)\in \tilde A}p_{r;i,i}^{(2n)}(x,x).$$
Consequently,
$$\sum_{k,\lambda_r(k)>0}\lambda_r(k)^{2n+1}\geq - \sum_{k,\lambda_r(k)<0}\lambda_r(k)^{2n+1}.$$
Since $|\lambda_r(k)|\leq 1$ for every $k$, we get
$$\sum_{k,\lambda_r(k)<0}\lambda_r(k)^{2n+2}\leq -\sum_{k,\lambda_r(k)<0}\lambda_r(k)^{2n+1}\leq \sum_{k,\lambda_r(k)>0}\lambda_r(k)^{2n+1},$$
hence
$$\sum_{k,\lambda_r(k)<0}\lambda_r(k)^{2n+2}\leq \sum_{k,\lambda_r(k)>0}\lambda_r(k)^{2n}.$$
We thus find
$$\mathrm{Tr}(Q_r^{2n+2})\leq 2 \sum_{k,\lambda_r(k)> 0}\lambda_r(k)^{2n}=2\sum_{k,\underset{\lambda(k)<1/2}{\lambda_r(k)>0,}}\lambda_r(k)^{2n}+2\sum_{k,\underset{\lambda(k)\geq 1/2}{\lambda_r(k)>0,}}\lambda_r(k)^{2n}.$$
Combining~(\ref{e:comparingeigenvalueslambdakgeq1/2}) and~(\ref{e:comparingeigenvalueslambdakleq1/2}), we finally deduce that
$$\mathrm{Tr}(Q_r^{2n+2})\leq 2 |\tilde A| (1-\epsilon(r))^{2n}+2\mathrm{Tr}(Q^{2\lfloor \epsilon(r)n\rfloor}).$$
This is true for any finite subset $A$ of $\Gamma$.
Therefore, Proposition~\ref{p:tracetheorem} shows that
$$\frac{1}{|E|}\sum_{i\in E}p_{r;i,i}^{(2n+2)}(e,e)\leq 2(1-\epsilon(r))^{2n}+2\frac{1}{|E|}\sum_{i\in E}\frac{1}{|E|}\mu_S^{(2\lfloor \epsilon(r)n\rfloor)}(e).$$
This concludes the proof of the first statement.
The second statement is obtained exactly the same way, using Lemma~\ref{l:epsilon1Dirichletnorms} instead of Lemma~\ref{l:epsilon0Dirichletnorms}.
\end{proof}

The original statement of Pittet and Saloff-Coste \cite{PSC} is as follows.
Given two quasi-isometric finitely generated amenable groups $\Gamma_1$ and $\Gamma_2$, endowed with simple random walks driven by probability measures $\mu_1$ and $\mu_2$, one has
$$\mu_1^{(2n+2)}(e)\lesssim (1-\epsilon)^{2n}+\mu_2^{(2\lfloor \epsilon n\rfloor )}(e).$$
Note that $\mu_1^{(2n)}$ and $\mu_2^{(2n)}$ decay sub-exponentially, as $\Gamma_1$ and $\Gamma_2$ are amenable.
Thus, one can get rid of the term $(1-\epsilon)^{2n}$ and so
$$\mu_1^{(2n+2)}(e)\lesssim \mu_2^{(2\lfloor \epsilon n\rfloor )}(e).$$
Similarly, one could get rid of the terms $(1-\epsilon(r))^{2n}$ in the statements of Proposition~\ref{p:continuousPSC}, since $p(r)$ is reversible and so by \cite[Corollary~12.12]{Woess-book}, its spectral radius is 1, hence $p_{r;i,i}^{(2n)}(e,e)$ also decays sub-exponentially.
However, in doing so, it is not clear whether the implicit constant would remain continuous.
In any case, combining Proposition~\ref{p:continuousPSC} and the original statement of Pittet and Saloff-Coste, we get the following result.
\begin{corollary}\label{c:continuousPSCdouble}
Let $\Gamma_1,\Gamma_2$ be two quasi-isometric finitely generated amenable groups and let $E_1,E_2$ be finite sets.
Let $r\mapsto p(r)$ be a continuous family of Markov transition kernels on $\Gamma_1\times E_1$ and $r\mapsto q(r)$ be a continuous family of Markov transition kernels on $\Gamma_2\times E_2$.
Assume that $p$ is $\Gamma_1$-invariant and $q$ is $\Gamma_2$-invariant and that both are strongly reversible, uniformly irreducible and have continuous second moment.
Then, there exists a continuous function $r\mapsto \epsilon(r)$ and there exists $n_0$ such that for every $n$
$$\frac{1}{|E_1|}\sum_{i\in E_1}p^{(2n)}_{r;i,i}(e,e)\lesssim \frac{1}{|E_2|}\sum_{i\in E_2}q^{(2\lfloor \epsilon(r)n\rfloor-2n_0)}_{r;i,i}(e,e)+(1-\epsilon(r))^{2n}$$
and
$$\frac{1}{|E_2|}\sum_{i\in E_2}q^{(2n)}_{r;i,i}(e,e)\lesssim \frac{1}{|E_1|}\sum_{i\in E_1}p^{(2\lfloor \epsilon (r)n\rfloor-2n_0)}_{r;i,i}(e,e)+(1-\epsilon(r))^{2n}.$$
\end{corollary}

\section{Applications to the first return kernel to a parabolic subgroup}\label{s:applicationsPSC}
Let $\Gamma$ be a relatively hyperbolic group.
In this section, we assume that parabolic subgroups are amenable.
Let $\mu$ be an admissible symmetric finitely supported probability measure on $\Gamma$.
As announced in the Introduction we assume further that $\mu(e)>0$ and we will explain at the end of the paper how to get rid of this assumption.
This will greatly simplify the results of this section.

We fix a parabolic subgroup $H$ and we fix $\eta\geq 0$.
The $\eta$-neighborhood $N_\eta(H)$ of $H$ can be identified with $H\times E$, where $E$ is a finite set.
Let $p(r)$ be the first return kernel to $N_\eta(H)$ associated with $r\mu$.
To simplify the notation, we will write, as before, $p_{r;i,j}$ instead of $p_{i,j}(r)$.

As explained, one of the major issues with the first return kernel is that it is not Markov.
The first step is to apply a suited transform.

\subsection{A Doob transform}
We consider the matrix $F(r)$ with entries
\begin{equation}\label{e:defmatrixF}
F_{i,j}(r)=\sum_{x\in H}p_{r;i,j}(e,x).
\end{equation}
This matrix was introduced in \cite{Martinfreeproducts}. It encodes many properties of the transition kernel $p(r)$.
It is called the vertical displacement matrix.
Let us recall some of its properties.

\begin{lemma}\label{l:Fstronglyirreducible}
The matrix $F(r)$ is strongly irreducible, i.e., there exists $n$ such that all its entries are positive.
\end{lemma}

\begin{proof}
In fact, the entries of $F(r)$ are given by
$$F_{i,j}^{(n)}(r)=p_{r;i,j}^{(n)}(e,x).$$
Since $p(r)$ is irreducible and $\mu(e)>0$, hence $p_{r;i,i}(e,e)>0$, there exists $n$ such that $p^{(n)}_{r;i,j}(e,e)>0$ for all $i,j$.
\end{proof}

According to the Perron-Frobenius theorem \cite[Theorem~1.1]{Seneta}, $F(r)$ has a dominant simple eigenvalue $\lambda(r)$ which is positive.
Letting $C(r)$ be a right-eigenvector, $C(r)$ has positive coefficients.
Moreover, any eigenvector of $F(r)$ with non-negative coefficients is in fact a multiple of $C(r)$ and is associated with the eigenvalue $\lambda(r)$.

The eigenvalue $\lambda(r)$ and the eigenvector $C(r)$ are closely related to the spectral data of $p(r)$.
In the context of virtually abelian parabolic subgroups, it is known that $\lambda(r)$ is, in fact, the spectral radius of $p(r)$.
We will prove below that this is still the case for arbitrary amenable parabolic subgroups.

\begin{lemma}\label{l:psymmetric}
The transition kernel $p(r)$ is symmetric and $H$-invariant.
\end{lemma}

\begin{proof}
Recall that given $x,y\in N_\eta(H)$,
$$p(r)(x,y)=G(x,y;N_\eta(H)^c|r),$$
that is,
$$
p(r)(x,y)=\ r\mu(x^{-1}y)+\sum_{n\geq 1}\sum_{z_1,...,z_n\notin N_\eta(H)}r^{n+1}\mu(x^{-1}z_1)\mu(z_1^{-1}z_2)...\mu(z_n^{-1}y).$$
Since $\mu$ is symmetric, any path $x,z_1,...,z_n,y$ avoiding $N_\eta(H)$ from $x$ to $y$ can be reversed to a path $y,z_n,...,z_1,x$ having the same weight 
$$\mu(y^{-1}z_n)\mu(z_n^{-1}z_{n-1})..\mu(z_{2}^{-1},z_1)\mu(z_1^{-1}x)=\mu(x^{-1}z_1)\mu(z_1^{-1}z_2)...\mu(z_{n-1}^{-1}z_{n})\mu(z_n^{-1}y).$$
Therefore,
$p(r)(x,y)=p(r)(y,x)$.
Similarly, $H$-invariance of $p(r)$ follows directly from $\Gamma$-invariance of $\mu$.
\end{proof}

\begin{remark}
Note that symmetry of $p(r)$ reads $p_{r;i,j}(x,y)=p_{r;j,i}(y,x)$.
In particular, $p_{r;i,j}$ is not itself symmetric.
\end{remark}

\begin{corollary}\label{c:Fsymmetric}
The matrix $F(r)$ is symmetric.
\end{corollary}

\begin{proof}
The entries of $F(r)$ are given by
$$F_{i,j}(r)=\sum_{x\in H}p_{r;i,j}(e,x).$$
Since $p(r)$ is symmetric and $H$-invariant, we get
$$F_{i,j}(r)=\sum_{x\in H}p_{r;i,j}(e,x)=\sum_{x\in H}p_{r;j,i}(x,e)=\sum_{x\in H}p_{r;j,i}(e,x^{-1})=\sum_{x\in H}p_{r;j,i}(e,x).$$
We thus find $F_{i,j}(r)=F_{j,i}(r)$.
\end{proof}

We denote by $c_i(r)$ the coordinates of $C(r)$, the eigenvector associated with $\lambda(r)$.
We now introduce the transition kernel $\tilde p(r)$ defined by
\begin{equation}\label{e:definitiontildep}
    \tilde p_{r;i,j}(x,y)=\frac{1}{\lambda(r) c_i(r)}p_{r;i,j}(x,y)c_j(r).
\end{equation}
In fact, $\tilde p$ is the Doob transform (also called $h$-process in literature) of $p$ with respect to the $\lambda(r)$-harmonic function $(x,i)\mapsto c_i(r)$.
We also introduce the function $m(r)$ on $N_\eta(H)=H\times E$ defined by
$$m_{r;i}(x):=m(r)((x,i))=c_i(r)^2.$$
Note that $m_{r;i}(x)$ is in fact independent of $x$.

\begin{lemma}\label{l:tildepMarkovinvariant}
The transition kernel $\tilde p(r)$ is Markov and $H$-invariant
\end{lemma}

\begin{proof}
Given $(x,i)$ in $H\times E$, we have
$$\sum_{y,j}\tilde p_{r;i,j}(r)(x,y)=\frac{1}{\lambda(r)c_i(r)}\sum_{y,j}p_{r;i,j}(x,y)c_j(r)$$
and by $H$-invariance of $p$,
$$\sum_{y,j}\tilde p_{r;i,j}(x,y)=\frac{1}{\lambda(r)c_i(r)}\sum_{y,j}p_{r;i,j}(e,x^{-1}y)c_j(r)=\frac{1}{\lambda(r)c_i(r)}\sum_{y,j}p_{r;i,j}(e,y)c_j(r).$$
Therefore,
$$\sum_{y,j}\tilde p_{r;i,j}(x,y)=\frac{1}{\lambda(r)}\frac{1}{c_i(r)}\sum_{j}F_{i,j}(r)c_j(r).$$
Since $C$ is an eigenvector for $F$ associated with $\lambda(r)$, we finally get
$$\sum_{y,j}\tilde p_{r;i,j}(x,y)=\frac{1}{\lambda(r)}\frac{1}{c_i(r)}\lambda(r)c_i(r)=1.$$
This proves that $\tilde p$ is Markov.
Now, $H$-invariance follows from $H$-invariance of $p$.
\end{proof}

\begin{lemma}\label{l:tildepreversible}
The transition kernel $\tilde p(r)$ is reversible with respect to the total conductance function $m(r)$.
\end{lemma}

\begin{proof}
Let $(x,i),(y,j)$ in $H\times E$.
Then,
$$m_{r;i}(x)\tilde p_{r;i,j}(x,y)=c_i(r)p_{r;i,j}(x,y)c_j(r)$$
and since $p(r)$ is symmetric,
$$m_{r;i}(x)\tilde p_{r;i,j}(x,y)=c_i(r)p_{r;j,i}(y,x)c_j(r),$$
hence
$m_{r;i}(x)\tilde p_{r;i,j}(x,y)=\tilde p_{r;j,i}(y,x)m_{r;j}(y)$.
\end{proof}

\subsection{The spectral radius and the vertical displacement matrix}

\begin{lemma}\label{l:sameconvolutionptildep}
The $n$th power of convolution $\tilde p^{(n)}$ of $\tilde p$ is given by $\widetilde{p^{(n)}}$, i.e.\ 
$$\tilde p^{(n)}_{r;i,j}(x,y)=\frac{1}{\lambda(r)^n c_i(r)}p^{(n)}_{r;i,j}(x,y)c_j(r).$$
\end{lemma}

\begin{proof}
We prove this by induction on $n$.
This is true by definition for $n=1$.
Let $n$ and assume that this is true for $n$.
Then,
$$\tilde p^{(n+1)}_{r;i,j}(x,y)=\sum_{z,k}\tilde p^{(n)}_{r;i,k}(x,z)\tilde p_{r;k,j}(z,y),$$
so
$$\tilde p^{(n+1)}_{r;i,j}(x,y)=\frac{1}{\lambda(r)^n}\frac{1}{c_i(r)}\sum_{z,k} p^{(n)}_{r;i,k}(x,y)c_k(r)\frac{1}{\lambda(r)c_k(r)}p_{r;k,j}(z,y)c_j(r).$$
Therefore,
$$\tilde p^{(n+1)}_{r;i,j}(x,y)=\frac{1}{\lambda(r)^{n+1}}\frac{1}{c_i(r)}p^{(n+1)}(x,y)c_j(r),$$
which concludes the proof.
\end{proof}

We let $\rho(r)$ be the spectral radius of $p(r)$.

\begin{proposition}\label{p:lambda=rho}
With the above notation,
$\rho(r)=\lambda(r)$.
\end{proposition}

\begin{proof}
Since $\tilde p(r)$ is a reversible Markov transition kernel on the amenable graph $N_\eta(H)=H\times E$, its own spectral radius is 1 by \cite[Corollary~12.12]{Woess-book}.
By Lemma~\ref{l:sameconvolutionptildep},
$$\tilde p^{(n)}_{r;i,i}(e,e)=\frac{1}{\lambda(r)^n c_i(r)}p^{(n)}_{r;i,i}(e,e)c_j(r).$$
Therefore,
$$\bigg( \tilde p^{(n)}_{r;i,j}(e,e)\bigg)^{1/n}\sim \frac{1}{\lambda(r)}\bigg(p^{(n)}_{r;i,i}(e,e)\bigg)^{1/n}.$$
Consequently, $1=\rho(r)/\lambda(r)$, that is $\rho(r)=\lambda(r)$.
\end{proof}

\subsection{Asymptotics of the powers of convolution}
We fix a finite generating set $S$ for $H$ and consider the associated simple random walk, which is driven by the uniform measure $\mu_H$ on $S$.

\begin{proposition}\label{p:tildepsatisfiesassumptionsPSC}
The map $r\in [1,R]\mapsto \tilde p(r)$ defines a continuous family of Markov transition kernels on $H\times E$ which is
\begin{itemize}
    \item $H$-invariant,
    \item strongly reversible with respect to the total conductance function $m(r)$,
    \item uniformly irreducible.
\end{itemize}
Moreover, if $\eta$ is large enough, it also has continuous second moment.
\end{proposition}

\begin{proof}
Recall that for $g,g'\in N_\eta(H)$,
$$p(r)(g,g')=G(g,g';N_\eta(H)^c|r).$$
Therefore, $p(r)$ is a power series with non-negative coefficients $a_n$ which are bounded by the coefficients $\mu^{(n)}$ of the Green function.
In particular, as a power series, $p(r)$ has a radius of convergence which is at least $R$ and so for every $g,g'\in N_\eta(H)$, $p(r)(g,g')$ is continuous in $r$.
Identifying $N_\eta(H)$ with $H\times E$, this means that for every $(x,i),(y,j)$, the function
$$r\mapsto p_{r;i,j}(x,y)$$
is continuous on $[1,R]$.
In other words, $r\mapsto p(r)$ is a continuous family of transition kernels.
We then deduce that the map $r\in [1,R]\mapsto F(r)$ is continuous.
Finally, since $F$ is strongly irreducible, it has a simple dominant eigenvalue.
By \cite[Theorem~VII.1.8]{Kato}, the maps $F\mapsto \lambda$ and $F\mapsto C$ are analytic; see also \cite[Lemma~3.3]{Martinfreeproducts}.
Consequently, the maps $r\mapsto \lambda(r)$ and $r\mapsto C(r)$ are continuous and so
$$r\mapsto \tilde p_{r;i,j}(x,y)$$
is continuous.
That is, $r\mapsto \tilde p(r)$ is a continuous family of Markov transition kernels.

Next, $H$-invariance follows from Lemma~\ref{l:tildepMarkovinvariant}.
Also, since $r\mapsto C(r)$ is continuous and $m(r)((x,i))$ depends only on $i$, the function $r\mapsto m_{r;i}(x)=1/c_i(r)$ is continuous and is uniformly bounded.
Thus, the fact that $r\mapsto \tilde p(r)$ is strongly reversible follows from Lemma~\ref{l:tildepreversible}.

We now prove that it is uniformly irreducible.
Recall that $S$ is a fixed finite generating set for $H$.
Since $\mu$ is admissible on $\Gamma$, for any $g,h$ in $\Gamma$ there exists $n_{g,h}$ such that $\mu^{(n_{g,h})}(g^{-1}h)$ is positive.
Assume that $g=(x,i)$ and $h=(y,j)$ are in $N_\eta(H)$.
Conditioning on successive visits to $N_\eta(H)$, we see that there exists $k_{g,h}$ such that $p^{(k_{g,h})}_{r;i,j}(x,y)$ is positive.
Also, since we assume in this section that $\mu(e)>0$, hence $p_{r;j,j}(y,y)>0$, we get that for any $k\geq k_{g,h}$, $p^{(k)}_{r;i,j}(x,y)$ is positive and so by Lemma~\ref{l:sameconvolutionptildep} $\tilde p^{(k)}_{r;i,j}(x,y)$ is positive.
We let $K$ be the maximum of the $k_{g,h}$, for $g=(e,i)$ and $h=(y,j)$ with $y\in S$ and we set
$$\epsilon(r)=\inf_{i,j}\inf_{y\in S}\tilde p ^{(K)}_{r;i,j}(e,y).$$
Then, $\epsilon$ is the infimum of a finite number of continuous functions, hence it is itself continuous.
This shows that $\tilde p$ is uniformly irreducible with respect to $\mu_H$.

Finally, we prove that for large enough $\eta$, $\tilde p$ has continuous second moment.
In fact, by \cite[Lemma~4.6]{DGstability}, for any $M$, if $\eta$ is large enough, then $p(r)$ has exponential moments up to $M$ for every $r$ and so has $\tilde p(r)$.
Therefore, the quantity $\tilde p_{r;i,j}(e,x)\mathrm{e}^{Md(e,x)}$ is uniformly bounded and so
$$\tilde p_{r;i,j}(e,x)\lesssim \mathrm{e}^{-Md(e,x)}$$

We now fix $M$ larger than the volume growth of $H$ and we fix $\eta$ large enough so that $p(r)$ has finite exponential moments up to $M$.
In particular, it has finite second moment and so $\tilde p$ also has finite second moment.
Now,
$$\mathcal{M}_2(\tilde p (r))=\sum_i\sum_{(x,j)\in N_\eta(H)}d(e,x)^2\tilde p_{r;i,j}(e,x).$$
Every term $d(e,x)^2\tilde p_{r;i,j}(e,x)$ is a continuous function and is bounded by what precedes by
$d(e,x)^2\mathrm{e}^{-Md(e,x)}$.
This quantity is summable, since $M$ is larger than the volume growth of $H$.
By dominated convergence, $\mathcal{M}_2(\tilde p (r))$ is continuous.
\end{proof}

We can now apply Proposition~\ref{p:continuousPSC} to the transition kernel $\tilde p(r)$ to prove the following.
In view of later applications, we use the notation $p_{r,H,\eta}$ instead of $p(r)$ for the first return transition kernel to $N_\eta(H)$, in order to insist on the involved parameters.
Similarly, we denote by $\rho_{H,\eta}(r)$ instead of $\rho(r)$ the spectral radius of $p_{r,H,\eta}$.
We also consider a finite generating set for $H$ and write $\mu_H$ for the uniform probability measure on $S$.
We are interested in asymptotics of $p_{r,H,\eta}^{(n)}$, as $n\to \infty$, for $r$ in a neighborhood of $R$.
Moreover, in order to get uniform bounds, we need to let $r$ vary in a compact set.
We arbitrarily fix the compact set $[1,R]$.

\begin{proposition}\label{p:applicationsonctinuousPSCfirstreturn}
Let $H$ be an amenable parabolic subgroup.
For large enough $\eta$, there exists a continuous function $r\mapsto \epsilon (r)\in (0,1/2)$, there exists $\kappa>0$ and there exists an integer $N$ such that for all $n$, for all $r\in [1,R]$,
$$\rho_{H,\eta}(r)^{-(2n+2)}p_{r,H,\eta}^{(2n+2)}(e,e)\lesssim \mu_H^{2\lfloor \epsilon(r)n\rfloor}(e)+\mathrm{e}^{-\kappa n}$$
and
$$\mu_H^{2\lfloor \epsilon(r)^{-1}n\rfloor+2}(e)\lesssim \rho_{H,\eta}(r)^{-2n-N}p_{r,H,\eta}^{(2n+N)}(e,e)+\mathrm{e}^{-\kappa n}.$$
\end{proposition}

\begin{proof}
We deduce from Proposition~\ref{p:continuousPSC} and Proposition~\ref{p:tildepsatisfiesassumptionsPSC} that
$$\sum_{i\in E}\tilde p^{(2n+2)}_{r;i,i}(e,e)\lesssim (1-\epsilon(r))^{2n}+ \mu_H^{2\lfloor \epsilon(r)n\rfloor}(e,e).$$
Identifying $N_\eta(H)$ with $H\times E$, we identify the neutral element $e$ of $\Gamma$ with $(e,i_0)$.
Recall that $C(r)$ has bounded positive coordinates and is continuous, hence $c_i(r)$ is bounded away from 0 and infinity on the compact set $[1,R]$.
By Lemma~\ref{l:sameconvolutionptildep}, we see that
$$\frac{1}{\lambda(r)^{2n+2}}p_{r,H,\eta}^{(2n+2)}(e,e)\lesssim \big(1-\epsilon(r)\big)^{2n}+\mu_H^{2\lfloor \epsilon(r)n\rfloor}(e,e).$$
Also, by continuity, $(1-\epsilon(r))^{2n}\leq \mathrm{e}^{-\kappa n}$ for some fixed $\kappa>0$.
Moreover, according to Proposition~\ref{p:lambda=rho}, $\lambda(r)=\rho_{H,\eta}(r)$.
Therefore,
$$\rho_{H,\eta}(r)^{-(2n+2)}p_{r,H,\eta}^{(2n+2)}(e,e)\lesssim \mu_H^{2\lfloor \epsilon(r)n\rfloor}(e,e)+\mathrm{e}^{-\kappa n}.$$
This proves the first inequality.
For the second one, we first use Proposition~\ref{p:continuousPSC} to get
$$\mu_H^{2\lfloor \epsilon(r)^{-1}n\rfloor+2}(e,e)\lesssim (1-\epsilon(r))^{2\lfloor \epsilon^{-1}(r) n\rfloor}+\sum_{i\in E}\tilde p_{r;i,i}^{(2n)}(e,e).$$
Since $\mu$ is admissible and $\mu(e)>0$, there exists $N$ such that for every $i$, for every $n$, $p_{r,H,\eta}^{(n)}((e,i),(e,i))\lesssim p_{r,H,\eta}^{(n+N)}((e,i_0),(e,i_0))$.
We thus get
$$\mu_H^{2\lfloor \epsilon(r)^{-1}n\rfloor+2}(e,e)\lesssim (1-\epsilon(r))^{2\lfloor \epsilon^{-1}(r)n\rfloor}+\rho_{H,\eta}(r)^{-2n-N}p_{r,H,\eta}^{(2n+N)}(e,e).$$
Finally, up to changing $\kappa$, we still have $(1-\epsilon(r))^{2\lfloor \epsilon^{-1}(r)n\rfloor}\leq \mathrm{e}^{-\kappa n}$.
This proves the second inequality.
\end{proof}

\begin{corollary}\label{c:LLTparabolicnilpotent}
Assume that parabolic subgroups are virtually nilpotent. Then for large enough $\eta$, for all $n$, for all $r\in [1,R]$,
$$p_{r,H,\eta}^{(n)}(e,e)\asymp \rho_{H,\eta}(r)^nn^{-d/2},$$
where the implicit constant is independent of $r$.
\end{corollary}

\begin{proof}
Consider a finite generating set for $H$ and let $\mu_H$ be the uniform measure on $S$.
Then, by Proposition~\ref{p:LLTnilpotent},
$$\mu_H^{*2n}(e)\asymp n^{-d/2}.$$
By Proposition~\ref{p:applicationsonctinuousPSCfirstreturn}, we see that
$$\rho_{H,\eta}(r)^{-(2n)}p_{r,H,\eta}^{(2n)}(e,e)\lesssim \mu_H^{2\lfloor \epsilon(r)(n-1)\rfloor}(e)+\mathrm{e}^{-\kappa n}\lesssim (\epsilon(r)(n-1))^{-d/2}+\mathrm{e}^{-\kappa n}.$$
Since $[1,R]$ is compact and $r\mapsto \epsilon(r)$ is continuous, it is bounded from above and below and so
$$\rho_{H,\eta}(r)^{-(2n)}p_{r,H,\eta}^{(2n)}(e,e)\lesssim n^{-d/2}+\mathrm{e}^{-\kappa n}$$
where the implicit constant is independent of $r$.
Also note that $\mathrm{e}^{-\kappa n}=o\big(n^{-d/2}\big)$ and so $n^{-d/2}+\mathrm{e}^{-\kappa n}\lesssim n^{-d/2}$.
Similarly,
$$(\epsilon(r)^{-1}(n-N))^{-d/2}\lesssim \rho_{H,\eta}(r)^{-(2n)}p_{r,H,\eta}^{(2n)}(e,e)+\mathrm{e}^{-\kappa n}$$ and so
$$n^{-d/2}\lesssim\rho_{H,\eta}(r)^{-(2n)}p_{r,H,\eta}^{(2n)}(e,e).$$
This proves the bounds for even $n$.
To conclude, note that $\alpha(r):=p_{r,H,\eta}(e,e)>0$ and that
$$\alpha(r) p_{r,H,\eta}^{(2n)}(e,e)\leq p_{r,H,\eta}^{(2n+1)}(e,e)\leq \frac{1}{\alpha(r)}p^{(2n+2)}_{r,H,\eta}(e,e).$$
Since $r\mapsto \alpha(r)$ is continuous, it is bounded from above and below and we get the desired bounds for odd $n$.
\end{proof}


\section{Asymptotics for the derivative of the spectral radius of the first return kernel}\label{s:derivativespectralradius}

Our next goal is to prove that the difference $G(e,e|R)-G(e,e|r)$ is roughly asymptotic to the difference $\rho_{H,\eta}(R)-\rho_{H,\eta}(r)$.
For that, we prove in fact that $\rho_{H,\eta}'(r)\asymp G'(e,e|r)$.
This will hold for every amenable parabolic subgroup $H$ and for every $\eta\geq 0$.
We now fix such a parabolic subgroup, we fix $\eta\geq 0$ and we write $p_r$ for the first return kernel to the $\eta$-neighborhood $N_\eta(H)$ and $\rho(r)$ for the associated spectral radius.
As above, we identify $N_\eta(H)$ with $H\times E$ where $E$ is a finite set and we write $p_{r;i,j}(x,y)=p_r((x,i),(y,j))$.
We let $P_r(x,y)$ be the matrix with coefficients $p_{r;i,j}(x,y)$.
We also let $G_r(x,y)$ be the matrix with coefficients $G_{r,H}((x,i),(y,j))$.

We will use the following set of inequalities, called relative Ancona inequalities.

\begin{theorem}\label{t:relativeAncona}
Let $\Gamma$ be a non-elementary relatively hyperbolic group.
Let $\mu$ be a symmetric, admissible, finitely supported probability measure on $\Gamma$.
Then, for every $c\geq 0$, there exists $C\geq 1$ such that the following holds.
For every $x,y,z$ such that $y$ is within $c$ of a relative geodesic from $x$ to $z$, for every $r\leq R$,
$$\frac{1}{C}G(x,y|r)G(y,z|r)\leq G(x,z|r)\leq C G(x,y|r)G(y,z|r).$$
\end{theorem}

This result has a rich history.
Those inequalities were first proved by Ancona \cite{Ancona} for $r=1$ for hyperbolic groups.
They are in fact the key result to prove many properties of random walks on hyperbolic groups.
They were extended to any $r\leq R$ with bounds not depending on $r$ by Gou\"ezel \cite{GouezelLLT} and Gou\"ezel-Lalley \cite{GouezelLalley}.
They were then proved for $r=1$ for relatively hyperbolic groups by Gekhtman-Gerasimov-Potyagailo-Yang \cite{GGPY} and finally extended to any $r\leq R$ in \cite{DGstability}.

We will especially apply relative Ancona inequalities in the following context.
Let $H$ be a parabolic subgroup and let $x\in \Gamma$.
Denote by $\pi_H(x)$ the nearest point projection of $x$ onto $H$.
Such a projection is well defined up to a bounded error, see \cite[Lemma~1.15]{Sisto-projections}.
Then combining  \cite[Lemma~1.13] {Sisto-projections} and  \cite[Lemma~1.15] {Sisto-projections}, any relatively geodesic from $e$ to $x$ passes within a bounded distance of $\pi_H(x)$.
Therefore, for every $r\leq R$,
\begin{equation}\label{e:Anconaprojectionparabolic}
G(e,x|r)\asymp G(e,\pi_H(x)|r)G(\pi_H(x),x|r),
\end{equation}
where the implicit constant does not depend on $r$.

\subsection{Rough asymptotics}
Recall that the vertical displacement matrix was introduced in Equation~(\ref{e:defmatrixF}).
It is strongly irreducible by Lemma~\ref{l:Fstronglyirreducible} and we denote by $\lambda(r)$ its dominant eigenvalue and by $C(r)$ an associated eigenvector.
We see $C(r)$ as a column vector and we write $C(r)^t$ for its transpose vector, which is thus a line vector.
We will write $A\cdot B$ for multiplication of matrices.
We normalize $C(r)$ to have $C(r)^t\cdot C(r)=1$.

We introduce the quantity
$$K(r)=\sum_n\sum_{x\in H}C(r)^t\cdot P_r^{(n)}(e,x)\cdot C(r).$$

\begin{lemma}\label{l:firstequationrhoandGreen}
We have
$$K(r)=\frac{1}{1-\rho(r)}=\sum_{x\in H}C(r)^t\cdot G_{r,H}(e,x)\cdot C(r).$$
\end{lemma}

\begin{proof}
For fixed $n$, we see that the matrix $\sum_xP_r^{(n)}(e,x)$ has coefficients
$$\sum_xp_{r;i,j}^{(n)}(e,x)=F_{i,j}^{(n)}(r).$$
In other words, $\sum_xP_r^{(n)}(e,x)=F(r)$.
Recall that $C(r)$ is both a right and left eigenvector for $F(r)$ associated with the eigenvalue $\lambda(r)$.
Also, recall that it satisfies $C(r)^t\cdot C(r)=1$, so we get
$$K(r)=\sum_n \lambda^n=\frac{1}{1-\lambda}.$$
Finally, since $\lambda=\rho$ by Proposition~\ref{p:lambda=rho}, this proves the first equality.
Now, since every term is non-negative in the double sum defining $K(r)$, we may invert the two sums, so that
$$K(r)=\sum_xC(r)^t\cdot \left (\sum_n P_r^{(n)}(e,x)\right)\cdot C(r)=\sum_xC(r) \cdot G_{r,H}(e,x)\cdot C(r),$$
proving the second equality.
\end{proof}

\begin{lemma}\label{l:differentiabilitysumGreen}
Given $i,j$, the function 
$$r\mapsto \sum_{x\in H} G_{r,H}((e,i),(x,j))$$
is differentiable on $(0,R)$ and its derivative is given by
$$\sum_{x\in H}G'_{r,H}((e,i),(x,j))$$
\end{lemma}

\begin{proof}
We write $G_{r;i,j}(e,x)$ for $G_{r,H}((e,i),(x,j))$.
According to Lemma~\ref{l:GreenderivativeI1},
$$G_{r;i,j}'(e,x)\asymp \sum_{y\in \Gamma}G(e,z|r)G(z,x|r).$$
Given $y\in H$, we let $\Gamma(y)$ be the set of $z\in \Gamma$ whose projection on $H$ is at $y$.
By relative Ancona inequalities~(\ref{e:Anconaprojectionparabolic}), for any $z\in \Gamma(y)$,
$$G(e,z|r)\asymp G(e,y|r)G(y,z|r)$$
and
$$G(z,x|r)\asymp G(z,y|r)G(y,e|r).$$
Therefore,
$$G_{r;i,j}'(e,x)\asymp \sum_{y\in H}G(e,y|r)G(y,x|r)\sum_{z\in \Gamma(y)}G(y,z|r)G(z,y|r).$$
By invariance by left translation, we see that
$$\sum_{z\in \Gamma(y)}G(y,z|r)G(z,y|r)\asymp \sum_{z\in \Gamma(e)}G(e,z|r)G(z,e|r).$$
Moreover, applying the above to $x=e$, we get
$$G'(e,e|r)\asymp  \sum_{z\in \Gamma(e)}G(e,z|r)G(z,e|r)I^{(1)}_H(r)$$
where we recall that
$$I^{(1)}_H(r)=\sum_{x\in H}G(e,x|r)G(x,e|r).$$
By Proposition~\ref{p:J1finite}, $I^{(1)}_H(r)$ is bounded from above and below.
This yields
$$G'(e,e|r)\asymp  \sum_{z\in \Gamma(e)}G(e,z|r)G(z,e|r)$$
and so
\begin{equation}\label{e:derivativeGreenparabolic}
G_{r;i,j}'(e,x)\asymp G'(e,e|r)\sum_{y\in H}G(e,y|r)G(y,x|r).
\end{equation}
Now fixing $r_1<R$ and considering a compact interval $[r_0,r_2]$ with $r_1\in (r_0,r_2)$ and $r_2<R$, we get for every $r\in [r_0,r_2]$
$$\sum_{x\in H}\sum_{y\in H}G(e,y|r)G(y,x|r)=\big(\sum_{x\in H}G(e,x|r)\big)^2\leq \big(\sum_{x\in H}G(e,x|r_2)\big)^2.$$
Note that the coordinates of $C(r)$ are bounded from above and below, so
$$\sum_{x\in H}G(e,x|r_2)\asymp C(r_2)^t\cdot \sum_{x\in H}G_{r_2,H}(e,x)\cdot C(r_2)$$
and so by Lemma~\ref{l:firstequationrhoandGreen}
$$\sum_{x\in H}\sum_{y\in H}G(e,y|r)G(y,x|r) \lesssim \frac{1}{(1-\rho(r_2))^2}.$$
By dominated convergence, we get the desired result at $r=r_1$ and since this is true for any $r_1<R$, this concludes the proof.
\end{proof}

\begin{proposition}\label{p:derivativeK}
We have
$$K'(r)=\frac{\rho'(r)}{(1-\rho(r))^2}=C(r)^t\cdot \left (\sum_{x\in H}G'_r(e,x)\right)\cdot C(r).$$
\end{proposition}

\begin{proof}
We saw in the proof of Proposition~\ref{p:tildepsatisfiesassumptionsPSC} that $r\mapsto C(r)$ is continuous.
In fact, this map is differentiable on $(0,R)$.
Indeed, by dominated convergence, $r\mapsto F(r)$ is differentiable on $(0,R)$ and once again, by \cite[Theorem~VII.1.8]{Kato}, the map $F\mapsto C$ is analytic.
Combining this with Lemma~\ref{l:differentiabilitysumGreen}, we see that
\begin{align*}
   K'(r)=&C'(r)^t\cdot \left (\sum_{x\in H}G_r(e,x)\right)\cdot C(r)+C(r)^t\cdot \left (\sum_{x\in H}G_r(e,x)\right)\cdot C'(r)\\
   &\hspace{3cm}+C(r)^t\cdot \left (\sum_{x\in H}G_r'(e,x)\right) \cdot C(r).
\end{align*}
Note that
$$\left (\sum_{x\in H}G_r(e,x)\right)\cdot C(r)=\sum_n\sum_x P^{(n)}_r(e,x)\cdot C(r).$$
We find that
$$\left (\sum_{x\in H}G_r(e,x)\right)\cdot C(r)= \sum_nF^n(r)\cdot C(r)=\frac{1}{1-\rho(r)} C(r).$$
Similarly,
$$C(r)^t\cdot \left (\sum_{x\in H}G_r(e,x)\right)=\frac{1}{1-\rho(r)}C(r)^t.$$
Therefore,
$$K'(r)=\frac{1}{1-\rho(r)}\bigg(C'(r)^t\cdot C(r)+C(r)^t\cdot C'(r)\bigg) +C(r)^t\cdot \left (\sum_{x\in H}G_r'(e,x)\right) \cdot C(r).$$
Finally, using that $C(r)^t\cdot C(r)=1$, we see that $C'(r)^t\cdot C(r)+C(r)^t\cdot C'(r)=0$ and so
$$K'(r)=C(r)^t\cdot \left (\sum_{x\in H}G'_r(e,x)\right)\cdot C(r).$$
This proves the second equality.
The first one follows immediately from the expression $K(r)=(1-\rho(r))^{-1}$.
\end{proof}

\begin{proposition}\label{p:rho'G'}
For every amenable parabolic subgroup $H$ and for every $\eta\geq 0$,
$$\rho_{H,\eta}'(r)\asymp G'(e,e|r).$$
\end{proposition}

\begin{proof}
By~(\ref{e:derivativeGreenparabolic}),
$$\sum_{x\in H}G_{r;i,j}'(e,x)\asymp G'(e,e|r)\sum_{x,y\in H}G(e,y|r)G(y,x|r).$$
Note that the coordinates of $C(r)$ are bounded from above and below, so
$$\sum_{x\in H}G(y,x|r)=\sum_{x\in H}G(e,x|r)\asymp C(r)^t\cdot G_r(e,x)\cdot C(r)=\frac{1}{1-\rho(r)}.$$
Therefore,
$$\sum_{x,y\in H}G(e,y|r)G(y,x|r)\asymp \frac{1}{(1-\rho(r))^2}.$$
We deduce from Proposition~\ref{p:derivativeK} that
$$\frac{\rho'(r)}{(1-\rho(r))^2}\asymp G'(e,e|r)\frac{1}{(1-\rho(r))^2}.$$
Consequently, $\rho'(r)\asymp G'(e,e|r)$.
\end{proof}


\begin{corollary}\label{c:Greenvsspectralradius}
For every amenable parabolic subgroup $H$ and for every $\eta\geq 0$,
$$G(R)-G(r)\asymp \rho_{H,\eta}(R)-\rho_{H,\eta}(r).$$
\end{corollary}










\subsection{Precise asymptotic for adapted random walks on free products}\label{ss:freeproductsrho'}
Let $H_0$ and $H_1$ be two finitely generated amenable groups.
We consider the free product $\Gamma=H_0*H_1$ and we set
$$\mu=\alpha \mu_0+(1-\alpha)\mu_1,$$
where $\mu_i$ is a finitely supported, admissible, symmetric probability measure on $H_i$.

Such a probability measure is called adapted to the free product structure.
Adapted probability measures on free products have been studied a lot and serve as a model for general random walks on relatively hyperbolic groups.
We refer to \cite[Chapter~17]{Woess-book} for more details.

As explained in the Introduction, free products are hyperbolic relative to the free factors $H_0$ and $H_1$.
We prove here the following.
\begin{proposition}\label{p:freeproductsrho'}
Let $\Gamma=H_0*H_1$ be a free product of finitely generated amenable groups.
Let $\mu$ be a finitely supported, symmetric and admissible adapted probability measure on a free product $\Gamma$.
Let $p_{H_i,r}$ be the first return kernel to $H_i$ and let $\rho_{H_i}$ be its spectral radius.
Then, there exists $C>0$ such that as $r$ tends to $R$,
$$\rho_{H_i}'(r)\sim CG'(e,e|r).$$
\end{proposition}

\begin{proof}
The main property of an adapted random walk on a free product is that it can only move inside one free factor at each step.
In particular, if the random walk leaves $H_0$ at some point $x$, it can only come back to $H_0$ at the same point.
Therefore, in this context, the first return kernel to $H_0$ is given by
$$p_{H_0,r}(e,\cdot)=w_0(r)\delta_e + r(1-\alpha)\mu_1(e)\delta_e+r\alpha \mu_0,$$
where $w_0$ is the weight of the first return to $e$ when taking the first step outside $H_0$.
Consequently, for $x\neq e$,
$$p_{H_0,r}'(e,x)=\alpha \mu_0(x)$$
and
$$p_{H_0,r}'(e,e)=w_0'(r)+(1-\alpha)\mu_1(e)+\alpha \mu_0(e).$$
Also, note that
$$\rho_{H_0}(r)=\sum_{x\in H_0}p_{H_0,r}(e,x).$$
Indeed, $p_{H_0,r}$ is a symmetric transition kernel on $H_0$ which is amenable.
Letting $t$ be its total mass, $t^{-1}p_{H_0,r}$ is Markov and symmetric and so its spectral radius is 1.
Therefore, the spectral radius of $p_{H_0,r}$ is its total mass $t$.
Thus, we get
\begin{equation}\label{e:comparingrhorandp'}
\rho'_{H_0}(r)=w_0'(r)+(1-\alpha)\mu_1(e)+\alpha=p_{H_0,r}'(e,e)+C.
\end{equation}

Now we want to compare $I^{(1)}(r)$ and $p_{H_0,r}'(r)$.
We start with the following equality:
\begin{equation}\label{equationGreenfirstreturne}
G(e,e|r)=1+\sum_{x\in H_0}p_{H_0,r}(e,x)G(x,e|r).
\end{equation}
This can be proved by conditioning on the first passage through $H_0$.
This can also be derived from
\begin{align*}
G(e,e|r)&=\sum_{n\geq 0}p_{H_0,r}^{(n)}(e,e)
=1+\sum_{n\geq 1}p_{H_0,r}*p_{H_0,r}^{(n-1)}(e,e)\\
&=1+\sum_{n\geq 0}\sum_{x\in H_0}p_{H_0,r}(e,x)p_{H_0,r}^{(n)}(x,e)
=1+\sum_{x\in H_0}p_{H_0,r}(e,x)\sum_{n\geq0}p_{H_0,r}^{(n)}(x,e)\\
&=1+\sum_{x\in H_0}p_{H_0,r}(e,x)G(x,e|r).
\end{align*}

Similarly, we have for any $x\neq e$ in $H_0$
\begin{equation}\label{equationGreenfirstreturnnote}
G(e,x|r)=\sum_{y\in H_0}p_{H_0,r}(e,y)G(y,x|r).
\end{equation}
According to Lemma~\ref{l:GreenderivativeI1},
multiplying~(\ref{equationGreenfirstreturne}) by $r$ and differentiating yields
$$I^{(1)}(r)=1+r\sum_{x\in H_0}p_{H_0,r}'(e,x)G(x,e|r)+\sum_{x\in H_0}p_{H_0,r}(e,x)I^{(1)}_{x^{-1}}(r).$$
By what precedes,
$$\sum_{x\in H_0}p_{H_0,r}'(e,x)G(x,e|r)=p_{H_0,r}'(e,e)G(e,e|r)+\alpha \sum_{x\neq e\in H_0}\mu_0(x)G(x,e|r).$$
Note that $\alpha\mu_0(x)\leq p_{H_0,r}(e,x)$ and so
$$\alpha\sum_{x\neq e\in H_0}\mu_0(x)G(x,e|r)\leq \sum_{x\in H_0}p_{H_0,r}(e,x)G(x,e|r)\leq G(e,e|r).$$
Therefore, the quantity
$$f(r)=\alpha\sum_{x\neq e\in H_0}\mu_0(x)G(x,e|r)$$
converges to some finite number $f(R)$ as $r$ tends to $R$.
We thus get
$$I^{(1)}(r)=rG(e,e|r)p_{H_0,r}'(e,e)+\sum_{x\in H_0}p_{H_0,r}(e,x)I^{(1)}_{x^{-1}}(r)+1+f(r).$$
Finally, by~(\ref{e:comparingrhorandp'})
\begin{equation}\label{equationI_1rho'}
I^{(1)}(r)=rG(e,e|r)\rho_{H_0}'(r)+\sum_{x\in H_0}p_{H_0,r}(e,x)I^{(1)}_{x^{-1}}(r)+1+f(r)-CrG(e,e|r).
\end{equation}

The next step is to find the exact correcting term when replacing $I^{(1)}_{x^{-1}}(r)$ by $I^{(1)}(r)$.
We have by definition
$$I^{(1)}_{x^{-1}}(r)=\sum_{y\in \Gamma}G(x,y|r)G(y,e|r).$$
Replacing the sum on $y$ by a sum on its projection $w$ on $H_0$ and a sum on $w^{-1}y=z$, we have
$$I^{(1)}_{x^{-1}}(r)=\sum_{w\in H_0}\sum_{z,\pi_{H_0}(z)=e}G(x,wz|r)G(wz,e|r).$$
Now, $G(x,wz|r)=\frac{G(x,w|r)G(e,z|r)}{G(e,e|r)}$ and
$G(wz,e|r)=\frac{G(z,e|r)G(w,e|r)}{G(e,e|r)}$.
Indeed, to get from $x$ to $wz$, the random walk needs to pass through the point $w$, as $\mu$ is adapted to the free product structure, so the random walk can only move inside one of the free factors at each step; see also \cite[Lemma~17.1~(b)]{Woess-book} for a formal proof.
This can be thought of as an exact version of the relative Ancona inequalities~(\ref{e:Anconaprojectionparabolic}) for adapted random walks on free products.
Consequently,
$$I^{(1)}_{x^{-1}}(r)=\sum_{z,\pi_{H_0}(z)=e}\frac{G(e,z|r)G(z,e|r)}{G(e,e|r)^2}\sum_{w\in H_0}G(x,w|r)G(w,e|r)$$
which can be written as
$$I^{(1)}_{x^{-1}}(r)=\sum_{z,\pi_{H_0}(z)=e}\frac{G(e,z|r)G(z,e|r)}{G(e,e|r)^2}I^{(1)}_{H_0,x^{-1}}(r).$$
This is true for any $x\in H_0$ and so in particular
$$I^{(1)}(r)=I^{(1)}_e(r)=\sum_{z,\pi_{H_0}(z)=e}\frac{G(e,z|r)G(z,e|r)}{G(e,e|r)^2}I^{(1)}_{H_0}(r),$$
hence
$$\frac{I^{(1)}_{x^{-1}}(r)}{I^{(1)}(r)}=\frac{I^{(1)}_{H_0,x^{-1}}(r)}{I^{(1)}_{H_0}(r)}.$$

Going back to~(\ref{equationI_1rho'}) and setting $g(r)=CG(e,e|r)-1-f(r)$ we find that
\begin{equation}\label{secondequationI_1rho'}
    \frac{\rho_{H_0}'(r)}{I^{(1)}(r)}rG(e,e|r)=1-\sum_{x\in H_0}p_{H_0,r}(e,x)\frac{I^{(1)}_{H_0,x^{-1}}(r)}{I^{(1)}_{H_0}(r)}+\frac{g(r)}{I^{(1)}(r)},
\end{equation}
where $g(r)$ converges to some $g(R)$ as $r$ tends to $R$.
To conclude, we estimate the quantity $\sum_{x\in H_0}p_{H_0,r}(e,x)I^{(1)}_{H_0,x^{-1}}(r)$.
This is by definition
$$\sum_{x\in H_0}\sum_{y\in H_0}p_{H_0,r}(e,x)G(x,y|r)G(y,e|r).$$
Inverting the two sums, we get
$$\sum_{y\in H_0}\bigg(\sum_{x\in H_0}p_{H_0,r}(e,x)G(x,y|r)\bigg) G(y,e|r).$$
Using~(\ref{equationGreenfirstreturne}) and~(\ref{equationGreenfirstreturnnote}), we see that this is in fact
$$\big(G(e,e|r)-1\big)G(e,e|r)+\sum_{y\neq e\in H_0}G(e,y|r)G(y,e|r)=I^{(1)}_{H_0}(r)-G(e,e|r).$$
To summarize,
$$\sum_{x\in H_0}p_{H_0,r}(e,x)I^{(1)}_{H_0,x^{-1}}(r)=I^{(1)}_{H_0}(r)-G(e,e|r).$$
Therefore,~(\ref{secondequationI_1rho'}) shows that
$$\frac{\rho_{H_0}'(r)}{I^{(1)}(r)}rG(e,e|r)=1-\frac{I^{(1)}_{H_0}(r)-G(e,e|r)}{I^{(1)}_{H_0}(r)}+\frac{g(r)}{I^{(1)}(r)},$$
that is,
$$\frac{\rho_{H_0}'(r)}{I^{(1)}(r)}rG(e,e|r)=\frac{G(e,e|r)}{I_1^{H_0}(e,e|r)}+\frac{g(r)}{I^{(1)}(r)}.$$
We conclude that
$\frac{\rho_{H_0}'(r)}{I^{(1)}(r)}$ converges to some finite quantity as $r$ tends to $R$, which concludes the proof.
For sake of completeness, let us compute the limit when the random walk is divergent.
Assuming further that $I^{(1)}(R)$ is infinite, we get
$$\frac{\rho_{H_0}'(r)}{I^{(1)}(r)}\underset{r\to R}{\longrightarrow}\frac{1}{RI^{(1)}_{H_0}(R)}.$$
Therefore,
$\rho_{H_0}'(r)\sim C_{H_0}I^{(1)}(r)$
where $C_{H_0}=\big (RI^{(1)}_{H_0}(R)\big)^{-1}$. 
\end{proof}

We conjecture that these precise asymptotics still hold in the general case.
Namely, given a non-elementary relatively hyperbolic group and a finitely supported admissible and symmetric probability measure on $\Gamma$, we conjecture that the following holds.
For any amenable parabolic subgroup $H$, for any $\eta\geq 0$, there exists $C$ such that as $r$ tends to $R$,
$$\rho_{H,\eta}'(r)\sim C G'(e,e|r).$$

\section{The convergent case}\label{s:LLTconvergent}
In this section we consider a relatively hyperbolic group $\Gamma$ with respect to virtually nilpotent subgroups.
We consider a finitely supported admissible and symmetric probability measure $\mu$ on $\Gamma$ and we assume that $I^{(1)}(R)$ is finite.

This implies in particular that the random walk is spectrally degenerate.
Therefore, $\mu$ has a well-defined rank of spectral degeneracy, which is the smallest possible homogeneous dimension of a parabolic subgroup $H$ along which the random walk is spectrally degenerate.
Our goal is to prove Theorem~\ref{t:maintheoremconvergent}, assuming further that $\mu(e)>0$.

\subsection{Weak Tauberian theorems}
We will use the following version of Karamata's Tauberian Theorem, which was proved in \cite{DussauleLLT1}, by using \cite[Theorem~2.10.2]{BinghamGoldieTeugels}.

\begin{proposition}\label{Karamata}\cite[Theorem~6.1]{DussauleLLT1}
Let $A(z)=\sum a_nz^n$ be a power series with non-negative coefficients $a_n$ and radius of convergence 1.
Let $\beta\geq 0$ and let $\ell$ be a slowly varying function.
Then, $\sum_{n\geq 0}a_ns^n  \asymp 1/(1-s)^\beta \ell(1/(1-s))$ for $s\in [0,1)$ if and only if $\sum_{k=0}^na_k\asymp \ell(n)n^{\beta}$.
The implicit constants are asked not to depend on $s$ and $n$ respectively.
%
\end{proposition}

We will also use the following technical lemma.
\begin{lemma}\label{l:sumvsleadingterm}
Let $q_n$ be a non-increasing sequence, let $\beta> 0$ and let $\ell$ be a slowly varying function.
If
$$\sum_{k=0}^nkq_k\asymp n^{\beta}\ell(n),$$
then
$$q_n\asymp n^{\beta-2}\ell(n).$$
\end{lemma}

\begin{proof}
Assume that
$$\sum_{k=0}^nkq_k\asymp n^{\beta}\ell(n).$$
Then, since $q_n$ is non-increasing,
$$\sum_{k=0}^nkq_n\leq \sum_{k=0}^{n}kq_k\lesssim n^{\beta}\ell(n).$$
Note that
$\sum_{k=0}^nk\asymp n^2$,
so
$n^2q_n\lesssim n^{\beta}\ell(n)$,
which proves the first inequality.
We now prove the second one.
Consider $c$ and $C$ such that
$$cn^{\beta}\ell(n)\leq \sum_{k=0}^{n}kq_k\leq Cn^{\beta}\ell(n).$$
fix a constant $A$ such that
$$4C\leq cA^{\beta}.$$
Since $\ell$ is slowly varying, $\ell(n)\leq 2\ell(An)$ for large enough $n$.
Therefore, if $n$ is large enough,
$$2Cn^{\beta}\ell(n)\leq c(An)^{\beta}\ell (An).$$
Consequently,
$$2Cn^{\beta}\ell(n)\leq \sum_{k=0}^{An}kq_k=\sum_{k=0}^nkq_k+\sum_{k=n+1}^{An}kq_k,$$
hence
$$2Cn^{\beta}\ell(n)\leq Cn^{\beta}\ell(n)+\sum_{k=n+1}^{An}kq_k.$$
Since $q_n$ is non-increasing,
$$\sum_{k=n+1}^{An}kq_k\leq  q_n\sum_{k=n+1}^{An}k\leq q_n C'n^2,$$
for some constant $C'$ which only depends on $A$.
This yields, for large enough $n$,
$$Cn^{\beta}\ell(n)\leq C'q_nn^2$$
proving the second inequality.
\end{proof}

Similarly, we have the following.
\begin{lemma}\label{l:sumvsleadingterm2}
Let $q_n$ be a non-increasing sequence, let $\beta> 0$ and let $\ell$ be a slowly varying function.
For every $j$, if
$$\sum_{k=0}^nk^jq_k\asymp n^{\beta}\ell(n),$$
then
$$q_n\asymp n^{\beta-j-1}\ell(n).$$
\end{lemma}


\subsection{Proof of the local limit theorem}
We now prove Theorem~\ref{t:maintheoremconvergent}, assuming that $\mu(e)>0$.
Recall that the quantities $I^{(k)}_x(r)$, $I^{(k)}_{H,\eta,x}(r)$ and $J^{(k)}(r)$ are defined in~(\ref{e:defIkg}),~(\ref{e:defIHk}) and~(\ref{e:defJeta}).
Also recall that $I^{(1)}(r)$ is the derivative of $rG(e,e|r)$.
As explained, similar formulae hold for higher derivatives. 
In fact, using \cite[Lemma~3.2]{DussauleLLT1}, $I^{(k)}(r)$ can be expressed as a sum:
\begin{equation}\label{e:IkGk}
I^{(k)}(r)=\sum_{i=0}^kP_i(r)G^{(i)}(e,e|r),
\end{equation}
where $P_i(r)$ is a polynomial function of $r$ of degree $i$.

\begin{lemma}
There exists $k$ such that $J^{(k)}(R)$ is infinite.
\end{lemma}

\begin{proof}
Let $H$ be a parabolic subgroup along which the random walk is spectrally degenerate.
Then, by Corollary~\ref{c:LLTparabolicnilpotent}, for large enough $\eta$, $p_{R,H,\eta}^{(n)}(e,e)\asymp n^{-d_H/2}$, where $d_H$ is the homogeneous dimension of $H$.
Consequently, setting $k=\lfloor d_H/2\rfloor$, the $k$th derivative of $G_{R,H,\eta}$ is infinite at 1, so combining  Lemma~\ref{l:GreenderivativeI1} and Lemma~\ref{l:sameGreen}, $I_{H,\eta}^{(k)}(R)$ is infinite, where we recall that
$$I_{H,\eta}^{(k)}(r)=\sum_{x_1,...,x_k\in N_\eta(H)}G(e,x_1|r)G(x_1,x_2|r)...G(x_{k-1},x_k|r)G(x_k,e|r).$$
Note that $I_{H,\eta}^{(k)}(r)\asymp I_H^{(k)}(r)$, hence $I_H^{(k)}(R)$ is infinite.
Thus, $J^{(k)}(R)$ is infinite.
\end{proof}

\begin{proposition}\label{p:IksameJk}
Let $k$ be the first integer such that $J^{(k)}(R)$ is infinite.
Then,
$$I^{(k)}(r)\asymp J^{(k)}(r).$$
\end{proposition}

\begin{proof}
First, for any $k$, $J^{(k)}(r)\leq I^{(k)}(r)$.
Second, assuming that $I^{(1)}(R)$ is finite, \cite[Lemma~5.7]{DussauleLLT1} shows that for every $k$, $I^{(k)}(r)$ is bounded by a sum only involving the quantities $I^{(j)}(r)$, $j<k$, and the quantities $J^{(j)}(r)$,  $j\leq k$.
Thus, by induction, $I^{(j)}(R)$ is finite for all $j<k$ and $I^{(k)}(r)\lesssim J^{(k)}(r)$.
\end{proof}

\begin{lemma}\label{l:IkJkd}
The first integer $k$ such that $J^{(k)}(R)$, equivalently $I^{(k)}(R)$, is infinite satisfies $d=2k+1$ or $d=2k+2$, where $d$ is the rank of spectral degeneracy of $\mu$.
\end{lemma}

\begin{proof}
Let $H$ be a parabolic subgroup.
Let $d_H$ be the homogeneous dimension of $H$ and let $k_H$ be the first integer such that $I_H^{(k)}(R)$ is infinite.
Then, by Corollary~\ref{c:LLTparabolicnilpotent}, $p_{R}^{(n)}(e,e)\asymp n^{-d_H/2}$.
First, assume that $d_H=2m+1$ is odd.
Then, for $l<m$, we have
$n^{l}p_{R}^{(n)}(e,e)\lesssim n^{l-m-1/2}\lesssim n^{-3/2}$, so $G^{(l)}_{H,R}(e,e)$ is finite.
Moreover, we have $n^{-1/2}\lesssim n^{m}p_R^{(n)}(e,e)$ and so $G^{(m)}_{H,R}(e,e)$ is infinite.
Therefore, by~(\ref{e:IkGk}), $k_H=m$.

Second, assume that $d_H=2m+2$ is even.
Similarly, for $l<m$, $G^{(l)}_{H,R}(e,e)$ is finite and $G^{(m)}_{H,R}(e,e)$ is infinite, so once again $k_H=m$.
We conclude that either $d_H=2k_H+1$ if $d_H$ is odd or $d_H=2k_H+2$ if $d_H$ is even.
This is true for every parabolic subgroup and so choosing the smallest possible $k_H$, we deduce that $d=2k+1$ if $d$ is odd and $d=2k+2$ if $d$ is even.
\end{proof}

We fix $\eta$ large enough, so that we can apply Corollary~\ref{c:LLTparabolicnilpotent}.
\begin{lemma}\label{l:estimateIkHoddandeven}
Let $H$ be a parabolic subgroup of homogeneous dimension $d_H$.
Write $d_H=2k+1$ if $d_H$ is odd and $d_H=2k+2$ if it is even.
Then,
$$I^{(k)}_H(r)\asymp \frac{1}{1-\rho_{H,\eta}(r)}$$
if $d_H$ is odd and
$$I^{(k)}_H(r)\asymp -\log (1-\rho_{H,\eta}(r))$$
otherwise.
\end{lemma}

\begin{proof}
By Corollary~\ref{c:LLTparabolicnilpotent}, we have
$$
p_{H,r,\eta}^{(n)}(e,e)\asymp \rho_{H,\eta}(r)^n n^{-d/2}.$$
Therefore,
$$n(n-1)...(n-k+1)p_{H,r,\eta}^{(n)}(e,e)\asymp n^kp_{H,r,\eta}^{(n)}(e,e)\asymp \rho_{H,\eta}(r)^n n^{-1/2}$$
if $d_H$ is odd and
$$n(n-1)...(n-k+1)p_{H,r,\eta}^{(n)}(e,e)\asymp n^kp_{H,r,\eta}^{(n)}(e,e)\asymp \rho_{H,\eta}(r)^n n^{-1}$$
if it is even.
In the former case, we deduce that
$$G^{(k)}(e,e|r)\asymp \sum_{n}\rho_{H,\eta}(r)^n n^{-1/2}\asymp \frac{1}{\sqrt{1-\rho_{H,\eta}(r)}}$$
and in the later case, we deduce that
$$G^{(k)}(e,e|r)\asymp \sum_{n}\rho_{H,\eta}(r)^n n^{-1}\asymp -\log (1-\rho_{H,\eta}(r).$$
Note that $k$ is the first integer such that $G^{(k)}(e,e|R)$ is infinite.
Consequently, $I^{(k)}(r)\asymp G^{(k)}(e,e|r)$ by~(\ref{e:IkGk}), which concludes the proof.
\end{proof}

We now let $d$ be the rank of spectral degeneracy of $\mu$ and we let $H_1,...,H_m$ be all the parabolic subgroups of homogeneous dimension $d$ along which $\mu$ is spectrally degenerate.
By what precedes,
\begin{equation}\label{e:IkJkforH1...Hm}
I^{(k)}(r)\asymp J^{(k)}(r)\asymp  \sum_{l=1}^m G^{(k)}_{H_l,r,\eta}(e,e|1).
\end{equation}
Also, by Corollary~\ref{c:LLTparabolicnilpotent}, we have
\begin{equation}\label{e:prhoforH1...Hm}
p_{H_l,r,\eta}^{(n)}(e,e)\asymp \rho_{H_l,\eta}(r)^n n^{-d/2}.
\end{equation}
To conclude the proof, we now treat separately the case where $d$ is odd and the case where $d$ is even.

\subsubsection{The odd case}
Assume that $d=2k+1$.
Then, by Lemma~\ref{l:estimateIkHoddandeven}
$$ G^{(k)}_{H_l,r,\eta}(e,e|1)\asymp \frac{1}{\sqrt{1-\rho_{H_l,\eta}(r)}}.$$
For $l=1,...,m$, $\mu$ is spectrally degenerate along $H_l$, so that $\rho_{H_l,\eta}(R)=1$, hence by Corollary~\ref{c:Greenvsspectralradius},
$$\sqrt{1-\rho_{H_l,\eta}(r)}\asymp \sqrt{G(e,e|R)-G(e,e|r)}.$$
Since $\mu$ is convergent, $G(e,e|R)-G(e,e|r)\asymp R-r$.
We thus get
$$\sqrt{1-\rho_{H_l,\eta}(r)}\asymp \sqrt{R-r}.$$
This is true for $l=1,...,m$ and so by~(\ref{e:IkJkforH1...Hm}),
$$I^{(k)}(r)\asymp \frac{1}{\sqrt{R-r}}.$$
Since $k$ is the first integer such that $I^{(k)}(R)$ is infinite, all derivatives $G^{(j)}(e,e|r)$ are finite at $R$ for $j<k$ and $I^{(k)}(r)\asymp G^{(k)}(e,e|r)$ by~(\ref{e:IkGk}).
We finally deduce that
\begin{equation}\label{e:convergentcaseasympGkdodd}
G^{(k)}(e,e|r)\asymp \frac{1}{\sqrt{R-r}}.
\end{equation}

We now apply Karamata's Tauberian Theorem.
We set $a_n=n^kR^n\mu^{(n)}(e)$.
For any $s<1$,
$$\sum_{n\geq 0}a_ns^n\asymp G^{(k)}(e,e|sR)\asymp \frac{1}{\sqrt{(1-s)}}.$$
By Proposition~\ref{Karamata}, we deduce that
$$\sum_{k=1}^na_k\asymp n^{1/2}.$$
Next, we want to apply Lemma~\ref{l:sumvsleadingterm2}, but
the sequence $a_k$ is not non-increasing.
However, by \cite[Corollary~9.4]{GouezelLalley}, using that $\mu$ is symmetric and aperiodic,
\begin{equation}\label{e:GouezelLalleydecreasing}
R^n\mu^{(n)}(e)=q_n+O\left (\mathrm{e}^{-\kappa n}\right )
\end{equation}
where $q_n$ is non-increasing and $\kappa>0$.
Since $R$ is the radius of convergence of the Green function, $R^n\mu^{(n)}(e)$ cannot decay exponentially fast and so
$$R^n\mu^{(n)}(e)\sim q_n.$$
We find that
$$\sum_{k=1}^nn^kq_k\asymp n^{1/2},$$
hence by Lemma~\ref{l:sumvsleadingterm2},
$$q_n\asymp n^{-k-1/2}=n^{-d/2}.$$
We use again that $R^n\mu^{(n)}(e)\asymp q_n$ to deduce that
$$\mu^{(n)}(e)\asymp R^{-n}n^{-d/2}.$$
This concludes the proof of Theorem~\ref{t:maintheoremconvergent} when $\mu(e)>0$ and the rank of spectral degeneracy $d$ is odd. \qed

\subsubsection{The even case}
Assume now that $d$ is even, so that $d=2k+2$.
Then by Lemma~\ref{l:estimateIkHoddandeven},
$$G^{(k)}_{H_l,r,\eta}(e,e|1)\asymp -\log (1-\rho_{H_l,\eta}(r))$$
and so we find this time
$$I^{(k)}(r)\asymp -\log (R-r).$$
Once again, $k$ is the first integer such that $I^{(k)}(R)$ is infinite and so by~(\ref{e:IkGk}),
$$G^{(k)}(e,e|r)\asymp I^{(k)}(r)\asymp -\log (R-r).$$

We would like to apply Karamata's Tauberian Theorem like in the odd case.
However, setting $a_k=n^kR^n\mu^{(n)}(e)$, we would find that $\sum_{1\leq k\leq n}a_k\asymp \log n$.
It is not clear if this would imply that $a_n\asymp 1/n$, even for a non-increasing sequence $a_k$.
To circumvent this issue, we go one step further and find the asymptotic behavior of $G^{(k+1)}(e,e|r)$.

Like in the proof of Proposition~\ref{p:IksameJk}, we use \cite[Lemma~5.7]{DussauleLLT1} to deduce that $I^{(k+1)}(r)$ is bounded by a sum of terms only involving the finite quantities $I^{(j)}(r)$ and $J^{(j)}(r)$, $j<k$, the quantities $I^{(k)}(r)J^{(2)}(r)$ and $I^{(2)}(r)J^{(k)}(r)$ and the quantity $J^{(k+1)}(r)$.
By~(\ref{e:IkJkforH1...Hm}), $I^{(k)}(r)$ is bounded by $J^{(k)}(r)$.
Similarly, $I^{(2)}(r)$ is bounded by $J^{(2)}(r)$. Indeed, either $k>2$ and so $I^{(2)}(R)$ and $J^{(2)}(R)$ are both finite, or $k=2$, since we must have $k>1$ by Proposition~\ref{p:J1finite}.
We thus get
$$I^{(k+1)}(r)\lesssim J^{(k+1)}(r)+J^{(k)}(r)J^{(2)}(r).$$
We now show that the dominant term is $J^{(k+1)}(r)$.
For the parabolic subgroups $H_1,...,H_m$, we have by~(\ref{e:prhoforH1...Hm})
$$n(n-1)...(n-k)p_{H_l,r,\eta}^{(n)}(e,e)\asymp n^{k+1}p_{H_l,r,\eta}^{(n)}(e,e)\asymp \rho_{H_l,\eta}(r)^n$$
and so
$$G^{(k+1)}_{H_l,r,\eta}(e,e|1)\asymp \frac{1}{1-\rho_{H_l,\eta}(r)}\asymp \frac{1}{R-r}.$$
For any other parabolic subgroup $H$, by what precedes, either $I_{H}^{(k+1)}(r)$ is finite, or $I_H^{(k+1)}(r)\asymp -\log (R-r)$, or $I_H^{(k+1)}(r)\asymp (R-r)^{-1/2}$.
Consequently,
$$J^{(k+1)}(r)\asymp \frac{1}{R-r}.$$
Next, we look at $J^{(k)}(r)J^{(2)}(r)$.
Either $k>2$ in which case $J^{(2)}(R)$ is finite and so
$J^{(k)}(r)J^{(2)}(r)\asymp -\log (R-r)$,
or $k=2$ in which case
$J^{(k)}(r)J^{(2)}(r)=J^{(k)}(r)^2$ and so
$J^{(k)}(r)J^{(2)}(r)\asymp \big(-\log (R-r)\big)^2$.
In both cases,
$J^{(k)}(r)J^{(2)}(r)\lesssim J^{(k+1)}(r)$.
Therefore, $I^{(k+1)}(r)\lesssim J^{(k+1)}(r)$.
Since by definition $J^{(k+1)}(r)\leq I^{(k+1)}(r)$, we deduce that
$$I^{(k+1)}(r)\asymp J^{(k+1)}(r)\asymp \frac{1}{R-r}.$$
This implies, in particular, that $G^{(k)}(e,e|r)=o\big(I^{(k+1)}(r)\big)$ and since $G^{(j)}(e,e|r)$ is finite for $j<k$, (\ref{e:IkGk}) shows that $I^{(k+1)}(r)\asymp G^{(k+1)}(e,e|r)$.
We finally get in the even case $d=2k+2$
\begin{equation}\label{e:convergentcaseasympGkdeven'}
G^{(k+1)}(e,e|r)\asymp \frac{1}{R-r}
\end{equation}

We can now apply Karamata's Tauberian Theorem to conclude.
We set this time $a_n=n^{k+1}R^n\mu^{(n)}(e)$.
For any $s<1$,
$$\sum_{n\geq 0}a_ns^n\asymp G^{(k)}(e,e|sR)\asymp \frac{1}{1-s} .$$
By Proposition~\ref{Karamata}, we deduce that
$\sum_{k=1}^na_k\asymp n$.

Once again, we cannot apply Lemma~\ref{l:sumvsleadingterm2} to the sequence $a_k$ which is not non-increasing, but we use~(\ref{e:GouezelLalleydecreasing}) to deduce that
$R^n\mu^{(n)}(e)\sim q_n$, where $q_n$ is non-increasing.
Then,
$$\sum_{k=1}^nn^kq_k\asymp  n$$
and by Lemma~\ref{l:sumvsleadingterm2},
$q_n\asymp n^{-k-1}=n^{-d/2}$.
We finally get
$$\mu^{(n)}(e)\asymp R^{-n}n^{-d/2}.$$
This concludes the proof of Theorem~\ref{t:maintheoremconvergent} when $\mu(e)>0$ and the rank of spectral degeneracy $d$ is even. \qed





\section{Divergent random walks with infinite Green moments}\label{s:LLTcritical}
In this section we finally assume that $I^{(1)}(R)$ is infinite and $J^{(2)}(R)$ is infinite.
That is, $\mu$ is divergent, but it is not spectrally positive recurrent.
We also assume that parabolic subgroups are amenable.
In Equation~(\ref{e:fundamental}), all three quantity diverge at $R$.
Therefore, both the underlying hyperbolic structure and the parabolic subgroups are influential in the asymptotic differential equation.

In what follows, whenever a parabolic subgroup $H$ is fixed, we will simply write $p_r$ for the first return kernel to $H$.

\subsection{Reduction to virtually nilpotent parabolic subgroups}\label{ss:reductionnilpotent}

\begin{lemma}\label{l:infinitederivativeimpliesnilpotent}
Let $H$ be a parabolic subgroup and let $k\geq 0$. If $I_{H}^{(k)}(R)$ is infinite, then $H$ is virtually nilpotent. Moreover, if $k$ is the smallest such integer, then $H$ has homogeneous dimension $2k+1$ or $2k+2$.
\end{lemma}

\begin{proof}
By a landmark result of Varopoulous \cite{Var86}, given a probability measure $\mu_0$ on a finitely generated group $\Gamma_0$, the following holds.
If $\Gamma_0$ is not virtually nilpotent, then for every $p>0$, $\mu_0^{(n)}(e)=O\big(n^{-p}\big)$.
We apply this to the first return kernel $p_R$ to $H$.
We assume that $I_{H}^{(k)}(R)$ is infinite for some $k$.
Then $G_{H,R}^{(k)}(e,e|1)$ is infinite.
This implies in particular that the spectral radius of $p_R$ is 1.
Recall that $p_R$ is $H$-invariant.
Denote by $t$ its total mass.
Then, $t^{-1}p_R$ is a Markov symmetric and $H$-invariant transition kernel on $H$.
Since $H$ is amenable, its own spectral radius is 1 and so $t=1$.
As a consequence, $p_R$ itself is Markov and thus defines a random walk on $H$.
By contradiction, if $H$ were not virtually nilpotent, then by the aforementioned Varopoulos' theorem it would hold that $n^kp_R^{(n)}(e,e)=O\big(n^{-2}\big)$ and so $G^{(k)}_{H,R}(e,e|1)$ would be finite.
We deduce that $H$ is indeed virtually nilpotent.
The second part of the statement now follows from Lemma~\ref{l:IkJkd}.
\end{proof}

The following corollary will allow us to reduce to the case where the parabolic subgroups are virtually nilpotent.

\begin{corollary}\label{c:J2infiniteimpliesnilpotency}
If $J^{(2)}(R)$ is infinite, then one of the parabolic subgroups along which $\mu$ is spectrally degenerate must be virtually nilpotent. Moreover, the rank of spectral degeneracy of $\mu$ must be 5 or 6.
\end{corollary}

\subsection{Proof of the local limit theorem}
We denote by $d$ the rank of spectral degeneracy of $\mu$, so we have $d=5$ or $d=6$.
We will treat these two cases separately.
In both cases, we first give an estimate for $J^{(2)}(r)$ and then we deduce from~(\ref{e:fundamental}) a new asymptotic differential equation satisfied by the Green function.
We then use Karamata's Tauberian Theorem to conclude.

In the remainder of the section, we fix $\eta$ large enough to apply Corollary~\ref{c:LLTparabolicnilpotent}.

\subsubsection{The even case $d=6$}
\begin{lemma}\label{l:asymptoticrank6}
Let $H$ be a virtually nilpotent parabolic subgroup of homogeneous dimension 6.
Assume that the random walk is spectrally degenerate along $H$.
Then, $I_{H}^{(2)}(r)\asymp - \log (G(R)-G(r))$.
\end{lemma}

\begin{proof}
For every $\eta \geq 0$ and $r\leq R$, we consider the first return kernel $p_{r,\eta}$ to $N_\eta(H)$.
Let $\rho_{\eta}(r)$ be the spectral radius of $p_{r,H}$ and let $R_\eta(r)$ be its inverse.
As a special case of Lemma~\ref{l:estimateIkHoddandeven}, for large enough $\eta$, we have
$I^{(2)}_H(r)\asymp -\log (1-\rho_{H,\eta}(r))$.
Since the random walk is spectrally degenerate along $H$, $\rho_{H,\eta}(R)=1$.
Corollary~\ref{c:Greenvsspectralradius} thus yields $I_{H,\eta}^{(2)}(r)\asymp -\log (G(R)-G(r))$.
For any given $\eta\geq 0$,
$I_{H,\eta}^{(2)}(r)\asymp I_{H}^{(2)}(r)$.
Consequently, we have
$I_{H}^{(2)}(r)\asymp -\log (G(R)-G(r))$.
\end{proof}

\begin{corollary}
Assume that $J^{(2)}(R)$ is infinite and
assume that the rank of spectral degeneracy of $\mu$ is 6.
Then, $J^{(2)}(r)\asymp -\log (G(R)-G(r))$.
\end{corollary}

\begin{proof}
Consider a parabolic subgroup $H$.
Then by Lemma~\ref{l:infinitederivativeimpliesnilpotent}, either $I_H^{(2)}(R)$ is finite or $H$ is virtually nilpotent of homogeneous dimension 5 or 6.
Since the possibility of $H$ having homogeneous dimension 5 is ruled out, we find that either $I_H^{(2)}(R)$ is finite or $I_H^{(2)}(r)\asymp -\log (G(R)-G(r))$ by Lemma~\ref{l:asymptoticrank6}.
Since there is at least one virtually nilpotent parabolic subgroup of homogeneous dimension 6, summing over all parabolic subgroups, we get
$J^{(2)}(r)\asymp -\log (G(R)-G(r))$.
\end{proof}

We can now complete the proof of the local limit theorem in the case $d=6$.
By the fundamental asymptotic differential equation~(\ref{e:fundamental}), we have
$$I^{(2)}(r)\asymp  I^{(1)}(r)^3 J^{(2)}(r)\asymp - I^{(1)}(r)^3 \log (G(R)-G(r)).$$
We rewrite this as
\begin{equation}\label{e:newdifferentialequationcriticalcasedeven}
\frac{I^{(2)}(r)}{I^{(1)}(r)^2}\asymp -I^{(1)}(r) \log (G(R)-G(r)).
\end{equation}
Recall that $I^{(1)}(r)\asymp G'(e,e|r)$ and that $I^{(2)}(r)\asymp G''(e,e|r)$.
Integrating~(\ref{e:newdifferentialequationcriticalcasedeven}) between $r$ and $R$, we get
$$\frac{1}{I^{(1)}(r)}\asymp - (G(R)-G(r))\log (G(R)-G(r)).$$
That is,
$$-I^{(1)}(r)(G(R)-G(r))\log (G(R)-G(r))\asymp 1$$
Integrating again between $r$ and $R$, we get
$$-\big(G(R)-G(r)\big)^2\log (G(R)-G(r))\asymp R-r.$$
Therefore,
$$-\big(G(R)-G(r)\big)\log (G(R)-G(r))\asymp \frac{R-r}{G(R)-G(r)}$$
and so
$$\frac{1}{I^{(1)}(r)}\asymp \frac{R-r}{G(R)-G(r)}.$$
We now rewrite this as
$$\frac{I^{(1)}(r) }{(G(R)-G(r))}\asymp \frac{1}{R-r}$$
and integrate once more between $r$ and $R$ to get
$$\log (G(R)-G(r))\asymp \log (R-r).$$
Finally,
\begin{align*}
\frac{1}{I^{(1)}(r)^2}&\asymp \bigg( -( G(R)-G(r))^2\log (G(R)-G(r))\bigg ) \bigg (-\log (G(R)-G(r))\bigg)\\
&\asymp -(R-r)\log (R-r),
\end{align*}
hence
\begin{equation}\label{e:criticalcaseasympG'even}
G'(e,e|r)\asymp \frac{1}{\sqrt{-(R-r)\log (R-r)}}.
\end{equation}

To conclude, we apply Karamata's Tauberian Theorem like in the convergent case.
We set $a_k=kR^k\mu^{(k)}(e)$.
Then, for any $s<1$,
$$\sum_{n\geq 0}a_ns^n\asymp I^{(1)}(sR)\asymp \frac{1}{\sqrt{-(1-s)\log (1-s)}}.$$
By Proposition~\ref{Karamata} applied to $a_k$ and to the slowly varying function $t\mapsto \log(t)^{-1/2}$ we thus have
$$\sum_{k=0}^n kR^k\mu^{(k)}(e)\asymp n^{1/2}\log (n)^{-1/2}.$$
By~(\ref{e:GouezelLalleydecreasing}),
$$R^n\mu^{(n)}(e)=q_n+O\left (\mathrm{e}^{-\kappa n}\right )$$
where $q_n$ is non-increasing and $\kappa>0$.
Since $R$ is the radius of convergence of the Green function $R^n\mu^{(n)}(e)$ cannot decay exponentially fast, so $R^n\mu^{(n)}(e)\sim q_n$ and
$$\sum_{k=0}^n k q_k\asymp n^{1/2}\log (n)^{-1/2}.$$
We thus get by Lemma~\ref{l:sumvsleadingterm}
$$q_n\asymp n^{-3/2}\log (n)^{-1/2}$$
and so
$$R^n\mu^{(n)}(e)\asymp n^{-3/2}\log (n)^{-1/2}.$$
We finally conclude that
$$\mu^{(n)}(e)\asymp R^{-n}n^{-3/2}\log (n)^{-1/2}.$$
This ends the proof of Theorem~\ref{t:maintheoremcritical} when $d$ is even and $\mu(e)>0$. \qed

\subsubsection{The odd case $d=5$}

\begin{lemma}\label{l:asymptoticrank5}
Let $H$ be a virtually nilpotent parabolic subgroup of homogeneous dimension 5.
Assume that the random walk is spectrally degenerate along $H$.
Then, $I_{H}^{(2)}(r)\asymp  \big (G(R)-G(r)\big)^{-1/2}$.
\end{lemma}

\begin{proof}
For every $\eta \geq 0$ and $r\leq R$, we consider the first return kernel $p_{r,\eta}$ to $N_\eta(H)$.
Let $\rho_{\eta}(r)$ be the spectral radius of $p_{r,H}$ and let $R_\eta(r)$ be its inverse.
As a special case of Lemma~\ref{l:estimateIkHoddandeven}, for large enough $\eta$, we have
$I^{(2)}(r)\asymp (1-\rho_{H,\eta}(r))^{-1/2}$.
Since the random walk is spectrally degenerate along $H$, $\rho_{H,\eta}(R)=1$.
Corollary~\ref{c:Greenvsspectralradius} thus yields
$I_{H,\eta}^{(2)}(r)\asymp \big (G(R)-G(r)\big)^{-1/2}$.
For any given $\eta\geq 0$,
$I_{H,\eta}^{(2)}(r)\asymp I_{H}^{(2)}(r)$.
Consequently, we have
$I_{H}^{(2)}(r)\asymp  \big (G(R)-G(r)\big)^{-1/2}$.
\end{proof}

\begin{corollary}
Assume that $J^{(2)}(R)$ is infinite and
assume that the rank of spectral degeneracy of $\mu$ is 5.
Then, $J^{(2)}(r)\asymp  \big (G(R)-G(r)\big)^{-1/2}$.
\end{corollary}

\begin{proof}
Consider a parabolic subgroup $H$.
Then by Lemma~\ref{l:infinitederivativeimpliesnilpotent}, either $I_H^{(2)}(R)$ is finite or $H$ is virtually nilpotent of homogeneous dimension 5 or 6.
We thus have three cases.
Either $I_H^{(2)}(R)$ is finite, or $I_H^{(2)}(r)\asymp -\log (G(R)-G(r))$ if the homogeneous dimension is 6 by Lemma~\ref{l:asymptoticrank6}, or
$I_H^{(2)}(r)\asymp \big (G(R)-G(r)\big)^{-1/2}$ if the homogeneous dimension is 5 by Lemma~\ref{l:asymptoticrank5}.
Note that
$$\log (G(R)-G(r))=O\bigg ( \big (G(R)-G(r)\big)^{-1/2}\bigg).$$
Since there is at least one virtually nilpotent parabolic subgroup of homogeneous dimension 5, we get
$J^{(2)}(r)\asymp  \big (G(R)-G(r)\big)^{-1/2}$ by summing over all parabolic subgroups.
\end{proof}

We can now complete the proof of the local limit theorem in the case $d=5$.
The proof is similar to that of case $d=6$, but simpler.
By the fundamental asymptotic differential equation~(\ref{e:fundamental}), we have
$$I^{(2)}(r)\asymp  I^{(1)}(r)^3 J^{(2)}(r)\asymp I^{(1)}(r)^3 \big (G(R)-G(r)\big)^{-1/2}.$$
We rewrite this as
\begin{equation}\label{e:newdifferentialequationcriticalcasedodd}
\frac{I^{(2)}(r)}{I^{(1)}(r)^2}\asymp I^{(1)}(r)\big (G(R)-G(r)\big)^{-1/2}.
\end{equation}
Recall that $I^{(1)}(r)\asymp G'(e,e|r)$ and that $I^{(2)}(r)\asymp G''(e,e|r)$.
Integrating~(\ref{e:newdifferentialequationcriticalcasedodd}) between $r$ and $R$, we find
$$\frac{1}{I^{(1)}(r)}\asymp \big (G(R)-G(r)\big)^{1/2}.$$
That is,
$$I^{(1)}(r)\big (G(R)-G(r)\big)^{1/2}\asymp 1.$$
Integrating again between $r$ and $R$, we get
$$\big (G(R)-G(r)\big)^{3/2}\asymp R-r,$$
hence
$$G(R)-G(r)\asymp (R-r)^{2/3}.$$
This yields
\begin{equation}\label{e:criticalcaseasympG'odd}
G'(e,e|r)\asymp (R-r)^{-1/3}.
\end{equation}

To conclude, we again apply Karamata's Tauberian Theorem.
As before, we set $a_k=kR^k\mu^{(k)}(e)$.
Then, for any $s<1$,
$$\sum_{n\geq 0}a_ns^n\asymp I^{(1)}(sR)\asymp (1-s)^{-1/3}.$$
By Proposition~\ref{Karamata} applied to $a_k$, we thus have
$$\sum_{k=0}^n kR^k\mu^{(k)}(e)\asymp n^{1/3}.$$
Once again, by \cite[Corollary~9.4]{GouezelLalley},
$$R^n\mu^{(n)}(e)=q_n+O\left (\mathrm{e}^{-\kappa n}\right )$$
where $q_n$ is non-increasing and $\kappa>0$.
Since $R$ is the radius of convergence of the Green function $R^n\mu^{(n)}(e)$ cannot decay exponentially fast, so $R^n\mu^{(n)}(e)\sim q_n$ and
$$\sum_{k=0}^n k q_k\asymp n^{1/3}.$$
We thus get by Lemma~\ref{l:sumvsleadingterm}
$$q_n\asymp n^{-5/3}$$
and so
$$R^n\mu^{(n)}(e)\asymp n^{-5/3}.$$
We finally conclude that
$$\mu^{(n)}(e)\asymp R^{-n}n^{-5/3}.$$
This ends the proof of Theorem~\ref{t:maintheoremcritical} when $d$ is odd and $\mu(e)>0$. \qed

\section{From lazy random walks to the general case}\label{s:lazynotlazy}
In Sections~\ref{s:LLTconvergent} and~\ref{s:LLTcritical}, we proved Theorem~\ref{t:maintheoremconvergent} and Theorem~\ref{t:maintheoremcritical} assuming further that the random walk is lazy, i.e., $\mu(e)>0$.
To conclude the proof of these theorems, we need first to treat the case of an aperiodic random walk which is not necessarily lazy.
Then, we need to treat the case of a random walk which is not even aperiodic.

\medskip
We first prove Theorem~\ref{t:maintheoremconvergent} and Theorem~\ref{t:maintheoremcritical} without assuming that $\mu(e)>0$.
Let $\mu$ be an admissible, symmetric, finitely supported probability measure and set
$$\tilde \mu(x)=\frac{1}{2}\mu(x)+\frac{1}{2}\delta_e(x).$$
Denote by $\tilde G$ the Green function associated with $\tilde \mu$.
Then,  by \cite[Lemma~9.2]{Woess-book},
$${\tilde{G}( x, y|r)=\frac{2}{2-r}G\bigg ( x, y\bigg|\frac{r}{2- r}\bigg ).}$$
We deduce that for every $k$, there exists a constant $C_k$ such that
\begin{equation}\label{e:GtildeGsamederivatives}
    G^{(k)}(r)\sim C_k\tilde G^{(k)}(r).
\end{equation}

The probability measure $\tilde \mu$ is also admissible, symmetric and finitely supported.
Moreover, it satisfies $\tilde \mu(e)>0$.
Using~(\ref{e:convergentcaseasympGkdodd}),~(\ref{e:convergentcaseasympGkdeven'}),~(\ref{e:criticalcaseasympG'even}) and~(\ref{e:criticalcaseasympG'odd}), we see that the following holds for $\tilde G$.
\begin{itemize}
    \item In the convergent case, letting $d$ be the rank of spectral degeneracy, writing $d=2k+2$ if $d$ is even and $d=2k+1$ otherwise, either
    $$\tilde G^{(k+1)}(e,e|r)\asymp (R-r)^{-1}$$
    if $d$ is even or
    $$\tilde G^{(k)}(e,e|r)\asymp (R-r)^{-1/2}$$
    if $d$ is odd.
    \item In the divergent, not spectrally positive recurrent case, letting $d$ be the rank of spectral degeneracy, either
    $$\tilde G'(e,e|r)\asymp \frac{1}{\sqrt{-(R-r)\log (R-r)}}$$
    if $d$ is even and
    $$\tilde G'(e,e|r)\asymp (R-r)^{-1/3}$$
    if $d$ is odd.
\end{itemize}

\medskip
Because of~(\ref{e:GtildeGsamederivatives}), the derivatives of the Green functions $G$ and $\tilde G$ have the same asymptotics.
In other words, the same results hold for the Green function $G$.

Therefore, in every situation encountered, we can apply once again Karamata's Tauberian Theorem to find the desired asymptotics for $\mu^{(n)}$.

\medskip
Also, let us explain how to treat the case of a random walk which is not aperiodic.
Note that $\mu$ is assumed to be symmetric.
Therefore, its period must be 1 or 2.
We see that even if $\mu$ is not aperiodic, $\mu^{(2)}$ is.
In fact, $\mu^{(2)}(e)\geq \mu(x)^2$ for any given $x$ in the support of $\mu$, so $\mu^{(2)}$ is lazy.
Therefore, we can apply the aperiodic case to $\mu^{(2)}$ to find the desired asymptotics for $\big(\mu^{(2)}\big)^{(n)}=\mu^{(2n)}$.

\section{The classification of local limit theorems and further remarks}\label{s:conclusionandclassification}
In this final section, we summarize the various existing results and we explain what is missing to obtain a complete classification of local limit theorems on relatively hyperbolic groups.

\subsection{Spectrally positive recurrent random walks}
We first recall the case of a spectrally positive recurrent random walk.
We assume that $I^{(1)}(R)$ is infinite and $J^{(2)}(R)$ is finite.

\begin{theorem}\cite[Theorem~1.4]{DussauleLLT1}
Let $\Gamma$ be a non-elementary relatively hyperbolic group and let $\mu$ be a finitely supported symmetric admissible probability measure on $\Gamma$.
Assume that $\mu$ is spectrally positive recurrent.
Then,
$$p_n(e,e)\asymp R^{-n}n^{-3/2}.$$
\end{theorem}

The proof of this theorem is very straightforward, once the fundamental asymptotic differential equation~(\ref{e:fundamental}) is established.
Indeed, assuming that $I^{(1)}(R)$ is infinite and $J^{(2)}(R)$ is finite, we have
$$I^{(2)}(r)\asymp I^{(1)}(r)^3.$$
By integrating between $r$ and $R$,
$$\frac{1}{I^{(1)}(r)^2}\asymp R-r,$$
hence $I^{(1)}(r)\asymp (R-r)^{-1/2}$ and so using Karamata's Tauberian Theorem, we find that $p_n(e,e)\asymp R^{-n}n^{-3/2}$.

Note that in the particular case of spectrally non-degenerate random walk, there is a more precise statement.
\begin{theorem}\cite[Theorem~1.1]{DussauleLLT2}
Let $\Gamma$ be a non-elementary relatively hyperbolic group and let $\mu$ be a finitely supported symmetric admissible probability measure on $\Gamma$.
Assume that $\mu$ is not spectrally degenerate.
Then, for every $x$, there exists $C_x$ such that as $n$ tends to infinity,
$$p_n(e,x)\sim C_x R^{-n}n^{-3/2}.$$
\end{theorem}

\subsection{Virtually abelian parabolic subgroups}
Whenever parabolic subgroups are in fact virtually abelian, the following holds for convergent random walks.

\begin{theorem}\cite[Theorem~1.3]{DPT}
Let $\Gamma$ be a non-elementary relatively hyperbolic group with respect to virtually abelian subgroups and let $\mu$ be a finitely supported symmetric admissible probability measure on $\Gamma$.
Assume that $\mu$ is convergent.
Let $d$ be the rank of spectral degeneracy of $\mu$.
Then, for every $x$ there exists $C_x$ such that as $n$ tends to infinity,
$$p_n(e,x)\sim C_x R^{-n}n^{-d/2}.$$
\end{theorem}

Let us now explain what is missing to prove precise asymptotics in the other cases.
First of all, in order to prove the local limit theorem with precise asymptotics in the spectrally positive recurrent case, one would need to prove that
$$I^{(2)}(r)\sim J^{(2)}(r)I^{(1)}(r)^3.$$
This in turn would prove that
$$\mu^{(n)}(e)\sim CR^{-n}n^{-3/2},$$
whether parabolic subgroups are virtually abelian or not.

In the critical case where $\mu$ is divergent but not spectrally positive recurrent, one really needs to find precise asymptotics of
the derivatives of the induced Green functions $t\mapsto G_{H,r,\eta}(e,e|t)$ at 1, as $r$ tends to $R$.
This requires precise asymptotics of the convolution powers of the first return kernel $p_{r,H,\eta}$.

It seems difficult to obtain such precise asymptotics if parabolic subgroups are not virtually abelian.
Indeed, recall that $p_{r,H,\eta}$ is not in general finitely supported but only has exponential moments (up to enlarging $\eta$).
Now, even for a fixed symmetric random walk with exponential moments on a nilpotent group, it is not proved that
$p^{(n)}(e,e)\sim Cn^{-d/2}$.

On the other hand, such asymptotics for the derivatives of the Green function $G_{H,r,\eta}$ have already been proved in \cite{DPT} when parabolic subgroups are virtutally abelian.
Therefore, the only missing step to prove the precise local limit theorem in the divergent, not spectrally positive recurrent case is that
$$\rho_{H,\eta}'(r)\sim C I^{(1)}(r),$$
which is the conjecture at the end of Section~\ref{ss:freeproductsrho'}.

\subsection{Free products}
Consider a free product
$\Gamma=H_0*H_1$ and an adapted probability measure
$$\mu=\alpha \mu_0+(1-\alpha)\mu_1.$$
In this situation, it is already known that
$$I^{(2)}(r)\sim I^{(1)}(r)^3J^{(2)}(r),$$
see for instance \cite[Proposition~4.1]{DPT23}.
Also, the convergent case is treated in \cite{CandelleroGilch}, not only for abelian free factors, but also for virtually nilpotent free factors.
Indeed, even though \cite[Theorem~1.1]{CandelleroGilch} only mentions groups of the form $\mathbb Z^{d_1}*\mathbb Z^{d_2}$, the authors treat the cases of probability measures for which the Green function has an expansion of the form
$$G(e,e|r)=a_0+a_1(R-r)+...+a_d(R-r)^d+a_{d+1}(R-r)^q\log^{\kappa}(R-r)+O\big((R-r)^{d+2}\big),$$
where $q\in (d,d+1]$ and $\kappa\in \mathbb Z_{\geq 0}$.
Since virtually nilpotent groups and abelian groups have the same type of local limit theorem with error terms, the Green functions present the same type of expansions and so the results of \cite{CandelleroGilch} also apply to virtually nilpotent free factors.

Finally,
according to Proposition~\ref{p:freeproductsrho'}, we do have
$$\rho_{H}'(r)\sim C I^{(1)}(r).$$
As a consequence, we can replace $\asymp$ with $\sim$ in Theorem~\ref{t:maintheoremcritical}, in the critical case of a divergent, not spectrally positive recurrent, random walk.

In summary, we have the following.
\begin{theorem}\label{t:classificationfreeproducts}
Let $H_0$ and $H_1$ be two virtually nilpotent finitely generated groups.
Let $\Gamma=H_0*H_1$ and let $\mu$ be an adapted probability measure on $\Gamma$.
Assume that $\mu$ is admissible, symmetric, finitely supported and aperiodic.
Then the following holds.
\begin{itemize}
\item If $\mu$ is convergent, then
    $$p_n(e,x)\sim C_x R^{-n}n^{-d/2},$$
    where $d$ is the rank of spectral degeneracy of $\mu$.
    \item If $\mu$ is divergent and spectrally positive recurrent, then
    $$p_n(e,x)\sim C_x R^{-n}n^{-3/2}.$$
    \item If $\mu$ is divergent and not spectrally positive recurrent, then its rank of spectral degeneracy $d$ must be 5 of 6.
    \begin{itemize}
        \item If $d=5$, then
        $$p_n(e,x)\sim C_x R^{-n}n^{-5/3}.$$
        \item If $d=6$, then
        $$p_n(e,x)\sim C_x R^{-n}n^{-3/2}\big(\log (n)\big)^{-1/2}.$$
    \end{itemize}
\end{itemize}
\end{theorem}

\subsection{Amenable parabolic subgroups}
We finally discuss the more general case where parabolic subgroups are only assumed to be amenable and we show limitation in our strategy.

First of all, in the spectrally positive recurrent case, the local limit theorem does not depend on parabolic subgroups and is always of the form
$\mu^{(n)}(e)\asymp R^{-n}n^{-3/2}$.
Second, Lemma~\ref{l:infinitederivativeimpliesnilpotent} has significant consequences. Indeed, as a particular case, it shows that if parabolic subgroups are not virtually nilpotent, the critical case of a divergent but not spectrally positive recurrent random walk cannot occur.

Therefore, assuming that parabolic subgroups are amenable but not virtually nilpotent, the only remaining case to treat is that of a convergent random walk.
Our application of Pittet-Saloff-Coste's Comparison Theorem in Proposition~\ref{p:applicationsonctinuousPSCfirstreturn} shows that the first return kernel to a parabolic subgroup behaves like a simple random walk, in terms of convolution powers.
It is reasonable to expect that $\mu^{(n)}R^{n}$ will also have the same behavior, since in the convergent case only the parabolic subgroups appear to be influential.
However, the strategy that we used for nilpotent parabolic subgroups cannot work in general.
Indeed, using again Lemma~\ref{l:infinitederivativeimpliesnilpotent}, we see that if parabolic subgroups are not virtually nilpotent, then every $J^{(k)}(R)$ is finite.
In particular, every $I^{(k)}(R)$ is finite and so all derivatives of the Green function are finite at the spectral radius.


\bibliographystyle{plain}
\bibliography{LLT_nilp}

\end{document}